\documentclass[a4paper,12pt]{article}     

\usepackage{graphicx}		  
\usepackage{amsmath,amssymb}	

\newtheorem{lemma}{Lemma}[section]
\newtheorem{theorem}[lemma]{Theorem}
\newtheorem{proposition}[lemma]{Proposition}
\newtheorem{corollary}[lemma]{Corollary}
\newtheorem{assumption}{Assumption}

\newcommand\LABEL[1]{\label{#1}}

\def\authorfont{\footnotesize}

\def\ccode#1{\par
\vspace*{8pt}
{\authorfont{\leftskip18pt\rightskip\leftskip
\noindent #1\par}}\par}

\newenvironment{Proof}{
\hspace*{-9mm}
{ \it Proof.}}
{\hfill {$\square$}\vspace{5mm}}

\begin{document}

\begin{center}{
{\Large 
 Properties of minimal charts and
 their applications VII: 
 charts of type $(2,3,2)$}
\vspace{10pt}
\\ 
Teruo NAGASE and Akiko SHIMA\footnote{The second author is supported by JSPS KAKENHI Grant Number 18K03309.}
}
\end{center}


\begin{abstract}
Let $\Gamma$ be a chart,
and we denote by $\Gamma_m$
the union of all the edges of label $m$.
A chart $\Gamma$ is of type $(2,3,2)$
if there exists a label $m$
such that 
$w(\Gamma)=7$,
$w(\Gamma_m\cap\Gamma_{m+1})=2$,
$w(\Gamma_{m+1}\cap\Gamma_{m+2})=3$,
and $w(\Gamma_{m+2}\cap\Gamma_{m+3})=2$
where 
$w(G)$ is the number of white vertices in $G$.
In this paper, we prove that there is 
no minimal chart of 
type $(2,3,2)$.
\end{abstract}

%
%
%
%

\ccode{2010 Mathematics Subject Classification. Primary 57Q45; Secondary 57Q35.}
\ccode{ {\it Key Words and Phrases}. surface link, chart, white vertex. }

\setcounter{section}{0}
\section{Introduction}


Charts are oriented labeled graphs in a disk (see  \cite{KnottedSurfaces},\cite{BraidBook}, and see Section~\ref{s:Prel}  for the precise definition of charts).
From a chart, we can construct an oriented closed surface 
embedded in 4-space ${\Bbb R}^4$ 
 (see \cite[Chapter 14, Chapter 18 and Chapter 23]{BraidBook}). 
A C-move 
is a local modification between two charts
in a disk (see Section~\ref{s:Prel} for C-moves).
A C-move between two charts induces 
an ambient isotopy between oriented closed surfaces 
corresponding to the two charts.

We will work in the PL category or smooth category. All submanifolds are assumed to be locally flat.
In \cite{ONS},
we showed that there is no minimal chart with exactly five vertices
 (see Section~\ref{s:Prel} for the precise definition of minimal charts). 
Hasegawa proved that there exists a minimal chart with exactly
six white vertices \cite{H1}. 
This chart represents a 2-twist spun trefoil.
In \cite{INS} and \cite{NST},
we investigated minimal charts with exactly four white vertices.
In this paper, 
we investigate properties of minimal charts 
which are used to prove that
there is no minimal chart with exactly seven white vertices
(see \cite{ChartApp1}, \cite{ChartAppII},
\cite{ChartAppIII}, \cite{ChartAppIV},
\cite{ChartAppV}, \cite{ChartAppVI},
\cite{ChartAppVIII}).

Let $\Gamma$ be a chart.
For each label $m$, we denote by $\Gamma_m$
the union of all the edges of label $m$.

Now we define a type of a chart:
Let $\Gamma$ be a chart with at least one white vertex, 
and $n_1,n_2,\dots,n_k$ integers.
The chart $\Gamma$ is of {\it type $(n_1,n_2,\dots,n_k)$} if there exists a label $m$ of $\Gamma$ satisfying the following three conditions:
\begin{enumerate}
\item[(i)] For each $i=1,2,\dots, k$, 
the chart $\Gamma$ contains exactly $n_{i}$ white vertices in $\Gamma_{m+i-1}\cap \Gamma_{m+i}$.
\item[(ii)] If $i<0$ or $i>k$, then $\Gamma_{m+i}$ does not contain any white vertices.
\item[(iii)] Both of the two subgraphs $\Gamma_m$ and $\Gamma_{m+k}$ contain at least one white vertex.
\end{enumerate}
If we want to emphasize the label $m$,
then we say that $\Gamma$ is of {\it type $(m;n_1,n_2,\dots,n_k)$}. 
Note that $n_1\ge1$ and $n_k\ge1$ by the condition (iii).

We proved in \cite[Theorem 1.1]{ChartAppII} that
if there exists a minimal $n$-chart $\Gamma$ with exactly seven white vertices,
then $\Gamma$ is a chart of 
type $(7),(5,2),(4,3),(3,2,2)$ or $(2,3,2)$ 
(if necessary we change the label
$i$ by $n-i$ for all label $i$).
In \cite{ChartAppV},
we showed that
there is no minimal chart of type $(3,2,2)$.
In \cite{ChartAppVI} and this paper,
we shall show 
the following.

\begin{theorem}
\LABEL{MainTheorem} 
There is 
no minimal chart of 
type $(2,3,2)$.
\end{theorem}

In the future paper \cite{ChartAppVIII},
we shall show there is no minimal chart of type
$(7),(5,2),(4,3)$.
Therefore we shall show that
there is no minimal chart with exactly seven white vertices.

The paper is organized as follows.
In Section~\ref{s:Prel},
we define charts and minimal charts.
In Section~\ref{s:ConnectedComponent},
we investigate the graphs $\Gamma_{m+1}$ and $\Gamma_{m+2}$
for a minimal chart $\Gamma$ of type $(m;2,3,2)$.
In Section~\ref{s:ROfamilies},
we consider a minimal chart $\Gamma$ of type $(m;2,3,2)$
such that $\Gamma_{m+1}$ contains either an oval
(a graph with two white vertices and two black vertices,
see Fig.~\ref{fig04}(a)), or 
one of the two graphs as shown in
Fig.~\ref{fig05}
(the graphs with five white vertices).
We investigate that the chart $\Gamma$ contains 
what kind of pseudo charts.
In Section~\ref{s:Ring},
we give two examples of non minimal charts $\Gamma$ of
type $(m;2,3,2)$
such that $\Gamma_{m+1}$ contains an oval.
In Section~\ref{s:IOC},
we review IO-Calculation
(a property of numbers of inward arcs of label $k$ and outward arcs of label $k$ in a closed domain $F$
with $\partial F\subset\Gamma_{k-1}\cup\Gamma_k\cup\Gamma_{k+1}$
for some label $k$).
In Section~\ref{s:XchangeLemma},
we investigate a minimal chart $\Gamma$
of type $(m;2,3,2)$
such that $\Gamma_{m+1}$ contains an oval.
We shall show that 
the chart $\Gamma$ contains the pseudo chart
as shown in Fig.~\ref{fig17}(a)
by C-moves without changing types of charts.
In Section~\ref{s:DiskLemma},
we shall show that
neither $\Gamma_{m+1}$ nor $\Gamma_{m+2}$
contains an oval 
for any minimal chart $\Gamma$ of type $(m;2,3,2)$.
In Section~\ref{s:3angledDisk},
we investigate a minimal chart $\Gamma$ 
of type $(m;2,3,2)$
such that $\Gamma_{m+1}$ contains 
the graph as shown in Fig.~\ref{fig05}(a).
In Section~\ref{s:CaseG},
we shall show that
neither $\Gamma_{m+1}$ nor $\Gamma_{m+2}$ contains
the graph as shown in Fig.~\ref{fig05}(a)
for any minimal chart $\Gamma$ of type $(m;2,3,2)$.
In Section~\ref{s:ShiftingLemma},
we give an example of a non minimal chart $\Gamma$
of type $(m;2,3,2)$
such that $\Gamma_{m+1}$ contains 
the graph as shown in Fig.~\ref{fig05}(b).
In Section~\ref{s:Main},
we shall prove Theorem~\ref{MainTheorem}.


\section{Preliminaries}
\LABEL{s:Prel}

In this section, 
we introduce 
the definition of charts and its related words.

Let $n$ be a positive integer.
An $n$-{\it chart}  
(a braid chart of degree $n$ \cite{KnottedSurfaces}
or a surface braid chart of degree $n$ \cite{BraidBook}) 
is 
an oriented labeled graph in the interior of a disk,
which may be empty 
or
have closed edges without vertices
satisfying the following four conditions
(see Fig.~\ref{fig01}):
\begin{enumerate}
\item[(i)] 
Every vertex has degree $1$, $4$, or $6$.
\item[(ii)] 
The labels of edges are 
in $\{1,2,\dots,n-1\}$.
\item[(iii)]
In a small neighborhood of
each vertex of degree $6$,
there are six short arcs,
three consecutive arcs are
oriented inward 
and
the other three are outward,
and
these six are labeled $i$ and $i+1$
alternately for some $i$,
where the orientation and label of
each arc are inherited from
the edge containing the arc.
\item[(iv)]
For each vertex of degree $4$,
diagonal edges have the same label
and
are oriented coherently,
and the labels $i$ and $j$ of
the diagonals satisfy $|i-j|>1$.
\end{enumerate}
We call a vertex of degree $1$ a {\it black vertex},
a vertex of degree $4$ a {\it crossing}, and 
a vertex of degree $6$ a {\it white vertex}
respectively.

Among six short arcs
in a small neighborhood of
a white vertex,
a central arc of each three consecutive arcs
oriented inward (resp. outward) 
is called a   
{\it middle arc} at the white vertex
(see Fig.~\ref{fig01}(c)).
For each white vertex $v$, 
there are two middle arcs at $v$ 
in a small neighborhood of $v$.
An edge is said to be {\it middle at} a white vertex $v$ if it contains a middle arc at $v$.

Let $e$ be an edge connecting $v_1$ and $v_2$.
If $e$ is oriented from $v_1$ to $v_2$,
then we say that 
$e$ is oriented {\it outward at $v_1$}
and {\it inward at $v_2$}


\begin{figure}[htb]
\begin{center}
\includegraphics{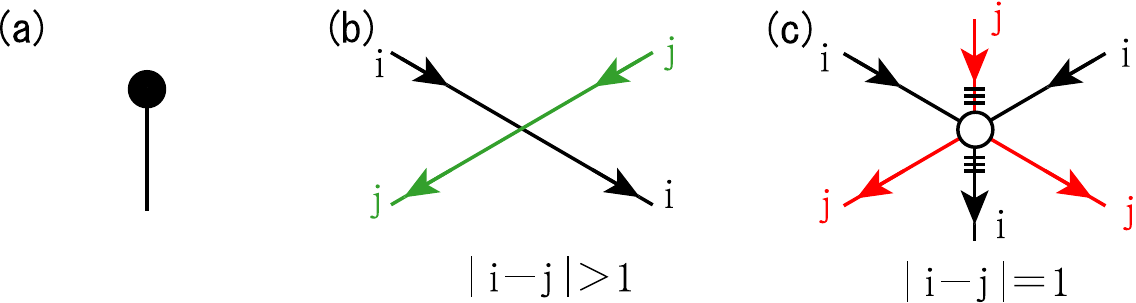}
\end{center}
\caption{ \LABEL{fig01} (a) A black vertex. (b) A crossing. (c) A white vertex. 
Each arc with three transversal short arcs is a middle arc at the white vertex. }
\end{figure}

Now {\it C-moves} are local modifications 
of charts as shown in Fig.~\ref{fig02}
(cf. \cite{KnottedSurfaces}, 
\cite{BraidBook} and \cite{Tanaka}).
Two charts are said to be {\it C-move equivalent}  if there exists
a finite sequence of C-moves 
which modifies one of the two charts 
to the other.

\begin{figure}
\begin{center}
\includegraphics{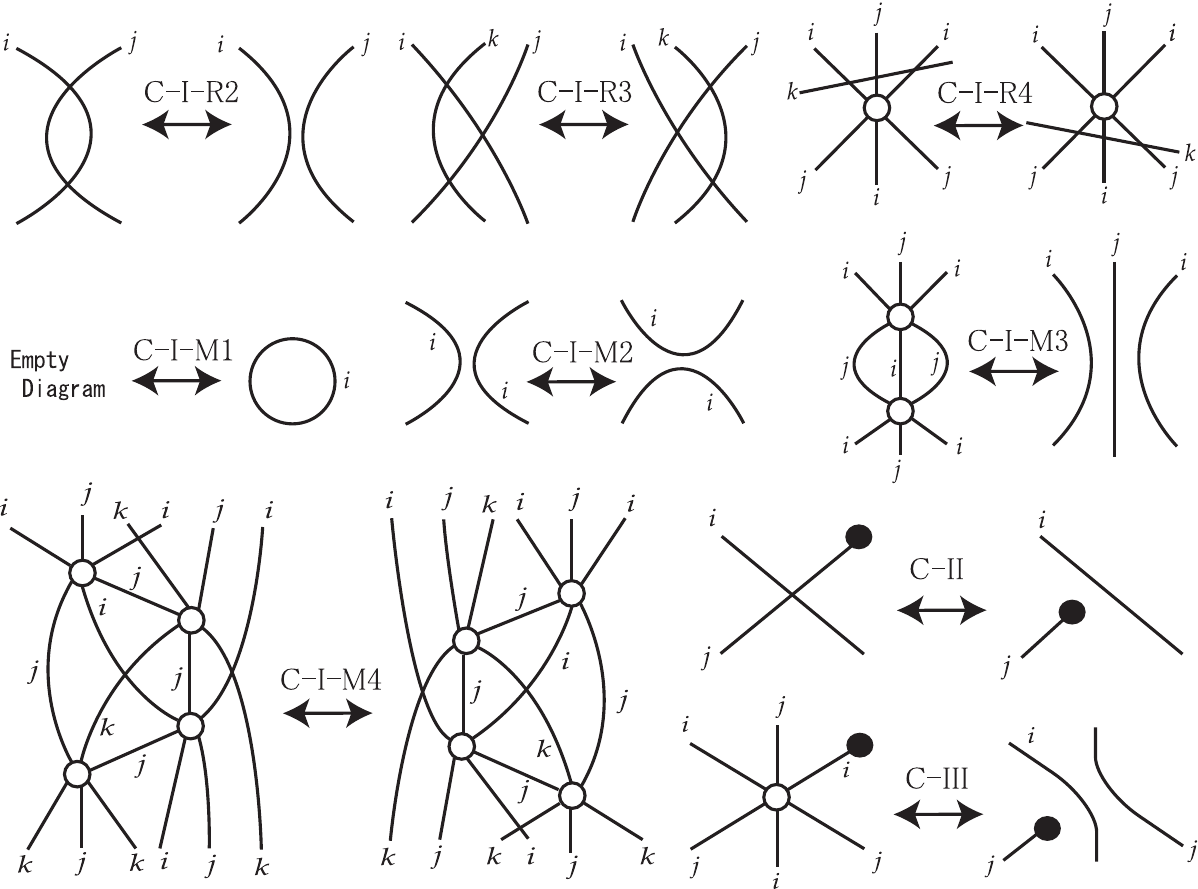}
\end{center}
\caption{ \LABEL{fig02} For the C-III move, 
the edge with the black vertex is not middle at
a white vertex in the left figure. }
\end{figure}

An edge in a chart is called 
a {\it free edge}
if it has
two black vertices.

For each chart $\Gamma$,
let $w(\Gamma)$ and $f(\Gamma)$ be the number of white vertices, and the number of free edges respectively.
The pair $(w(\Gamma), -f(\Gamma))$ is called a {\it complexity} of the chart (see \cite{BraidThree}).
A chart $\Gamma$ is called a {\it minimal chart} if its complexity is minimal among the charts C-move equivalent to the chart $\Gamma$ with respect to the lexicographic order of pairs of integers.

We showed the difference of a chart in a disk and in a 2-sphere (see \cite[Lemma 2.1]{ChartApp1}).
This lemma follows from that there exists a natural one-to-one correspondence between $\{$charts in $S^2\}/$C-moves and $\{$charts in $D^2\}/$C-moves, conjugations
(\cite[Chapter 23 and Chapter 25]{BraidBook}).
To make the argument simple, we assume that 
the charts lie on the 2-sphere instead of the disk.
\begin{assumption}
In this paper,
all charts are contained in the $2$-sphere $S^2$.
\end{assumption}
We have the special point in the 2-sphere $S^2$, called the point at infinity,
 denoted by $\infty$.
In this paper, all charts are contained in a disk such that the disk 
does not contain the point at infinity $\infty$.

An edge in a chart is called 
a {\it terminal edge}
if it has
a white vertex and a black vertex.

Let $\Gamma$ be a chart,
and $m$ a label of $\Gamma$. 
A {\it hoop} is a closed edge of $\Gamma$ without vertices 
(hence without crossings, neither).
A {\it ring} is a simple closed curve in $\Gamma_m$ containing a crossing but not containing any white vertices.
A hoop is said to be {\it simple} 
if one of the two complementary domains
of the hoop
does not contain any white vertices.

We can assume that
all minimal charts $\Gamma$
satisfy the following four conditions 
(see \cite{ChartApp1},\cite{ChartAppII},\cite{ChartAppIII},
\cite{StI}):

\begin{assumption}
\LABEL{AssumeTerminal}
If an edge of $\Gamma$
contains a black vertex,
then the edge is a free edge 
or a terminal edge.
Moreover 
any terminal edge contains a middle arc.
\end{assumption}

\begin{assumption}
\LABEL{NoSimpleHoop}
All free edges and simple hoops in $\Gamma$ 
are moved into a small neighborhood $U_\infty$ 
of the point at infinity $\infty$. 
Hence
we assume that 
$\Gamma$ does not contain free edges
nor simple hoops, 
otherwise mentioned. 
\end{assumption}

\begin{assumption}
\LABEL{Ring}
Each complementary domain of
any ring and hoop must contain 
at least one white vertex. 
\end{assumption}

\begin{assumption}
\LABEL{Infinity}
The point at infinity $\infty$ is moved in any complementary domain of $\Gamma$.
\end{assumption}

In this paper
for a set $X$ in a space
we denote 
the interior of $X$,
the boundary of $X$ and
the closure of $X$
by Int$X$, $\partial X$
and $Cl(X)$
respectively.

Let $\alpha$ be a simple arc or an edge,
and $p,q$ the endpoints of $\alpha$.
We denote 
$\partial \alpha=\{p,q\}$ and ${\rm Int}\alpha=\alpha-\partial \alpha$.

\section{The graphs $\Gamma_{m+1}$ and $\Gamma_{m+2}$}
\LABEL{s:ConnectedComponent}
In this section,
we investigate the graphs $\Gamma_{m+1}$ and $\Gamma_{m+2}$
for a minimal chart $\Gamma$ of type $(m;2,3,2)$.

 In our argument  we often construct a chart $\Gamma$. 
On the construction of a chart $\Gamma$, for a white vertex $w\in\Gamma_m$ for some label $m$,  
among the three edges of $\Gamma_m$ 
containing $w$, 
if one of the three edges is a terminal edge 
(see Fig.~\ref{fig03}(a) and (b)), 
then we remove the terminal edge and
put a black dot at the center of the white vertex  as shown in Fig.~\ref{fig03}(c).
Namely
Fig.~\ref{fig03}(c) means 
Fig.~\ref{fig03}(a) or 
Fig.~\ref{fig03}(b).
We call the vertex in Fig.~\ref{fig03}(c) 
a {\it BW-vertex} with respect to $\Gamma_m$.

\begin{figure}[htb]
\centerline{\includegraphics{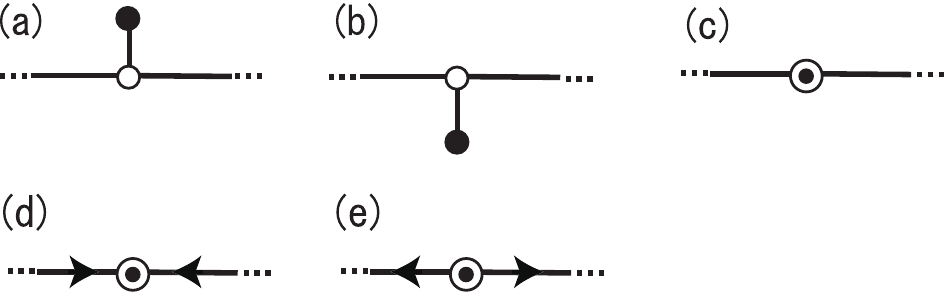}}
\caption{\LABEL{fig03}
(a),(b) White vertices in terminal edges.
(c),(d),(e) BW-vertices.}
\end{figure}

\begin{lemma}
\LABEL{OriBWvertex}
{\rm (\cite[Lemma 3.1]{ChartAppV})}
In a minimal chart $\Gamma$,
for each BW-vertex in $\Gamma_m$,
the two edges of label $m$ containing the BW-vertex
are oriented inward or outward at the BW-vertex
simultaneously
if each of the two edges is not a terminal edge
$($see Fig.~\ref{fig03}$($d$)$ and $($e$))$.
\end{lemma}

Let $X$ be a set in a chart $\Gamma$.
Let
 $$w(X)=\text{the number of white vertices in $X$.}$$

\begin{lemma}
\LABEL{LemmaWithTerminal3}
{\rm (\cite[Lemma 3.2(1)]{ChartAppV})}
Let $\Gamma$ be a minimal chart,
and $m$ a label of $\Gamma$.
Let $G$ be a connected component of $\Gamma_m$.
If $1\le w(G)$, then $2\le w(G)$.
\end{lemma}

We call 
an {\it oval},
a {\it skew $\theta$-curve} the two graphs as shown
in Fig.~\ref{fig04}(a),(b)
respectively.

\begin{figure}
\centerline{\includegraphics{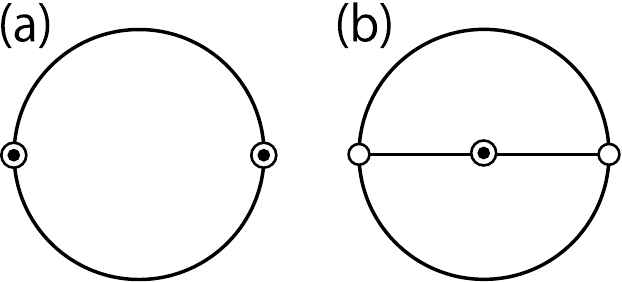}}
\caption{\LABEL{fig04}
(a) An oval. (b) A skew $\theta$-curve.}
\end{figure}

\begin{lemma}
{\rm (\cite[Theorem 1.3]{ChartAppVI})}
\LABEL{OvalTypeGTypeH}
If there exists a minimal chart $\Gamma$ of 
type $(m;2,3,2)$,
then each of $\Gamma_{m+1}$ and $\Gamma_{m+2}$
contains either the union of 
an oval and a skew $\theta$-curve, or
 one of two graphs as shown 
in  Fig.~\ref{fig05}.
\end{lemma}

\begin{figure}
\centerline{\includegraphics{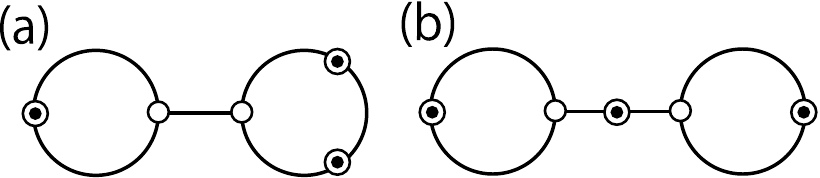}}
\caption{\LABEL{fig05}
Graphs with three black vertices.
}
\end{figure}


\section{RO-families of pseudo charts}
\LABEL{s:ROfamilies}

In this section
we consider a minimal chart $\Gamma$ of type $(m;2,3,2)$
such that $\Gamma_{m+1}$ contains either an oval, or 
one of the two graphs as shown in
Fig.~\ref{fig05}.
We investigate that the chart $\Gamma$ contains 
what kind of pseudo charts.

Let $\Gamma$ be a chart. 
Suppose that an object consists of 
some edges of $\Gamma$, arcs in edges of $\Gamma$ and arcs around white vertices.
Then the object is called a {\it pseudo chart}.

Let $\Gamma$ be a chart, 
$D$ a disk, and 
$G$ a pseudo chart with $G \subset D$.
Let $r:D\to D$ be a reflection of $D$, and $G^*$ the pseudo chart obtained from $G$ by changing the orientations of all of the edges.
Then the set $\{G,G^*, r(G), r(G^*)\}$ 
is called the {\it RO-family of the pseudo chart $G$}.

In our argument,
we often need a name for an unnamed edge by using a given edge and a given white vertex.
For the convenience,
we use the following naming:
Let $e',e_i,e''$ be three consecutive edges containing  a white vertex $w_j$. Here, 
the two edges $e'$ and $e''$ are unnamed edges. 
There are six arcs in a neighborhood $U$ of the white vertex $w_j$. 
If the three arcs $e'\cap U$, $e_i \cap U$, $e'' \cap U$ lie anticlockwise around the white vertex $w_j$ in this order, 
then $e'$ and $e''$ are denoted by $a_{ij}$ and $b_{ij}$ 
respectively (see Fig.~\ref{fig06}).
There is a possibility $a_{ij}=b_{ij}$ if they contain only one white vertex.

\begin{figure}[htb]
\centerline{\includegraphics{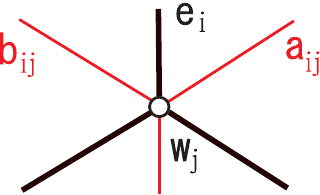}}
\caption{\LABEL{fig06}
The three edges $a_{ij},e_i,b_{ij}$ are consecutive edges around the white vertex $w_j$.}
\end{figure}

\begin{lemma}
{\rm (\cite[Lemma 8.3]{ChartAppVI})}
\LABEL{OvalGammaM+1Step0}
Let $\Gamma$ be a minimal chart 
of type $(m;2,3,2)$.
If $\Gamma_{m+1}$ contain an oval,
then 
$\Gamma$ contains one of the RO-families of
the two pseudo charts
as shown in Fig.~\ref{fig07}.
\end{lemma}

\begin{figure}
\centerline{\includegraphics{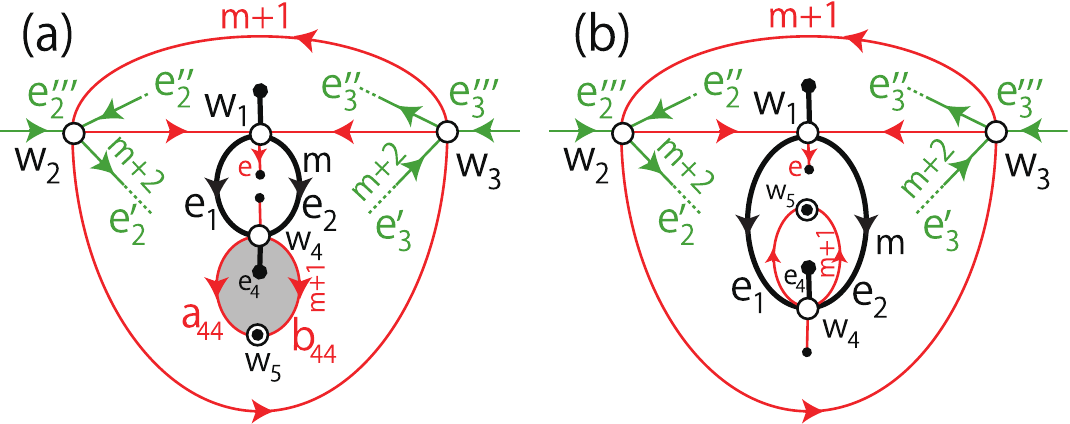}}
\caption{
\LABEL{fig07}
Pseudo charts containing
a skew $\theta$-curve and an oval of label $m+1$.
}
\end{figure}

\begin{lemma}
{\rm (\cite[Lemma 8.4]{ChartAppVI})}
\LABEL{LemmaTypeG}
Let $\Gamma$ be a minimal chart 
of type $(m;2,3,2)$.
If $\Gamma_{m+1}$ 
contains the graph as shown in 
Fig.~\ref{fig05}$($a$)$,
then 
$\Gamma$ contains one of the RO-families of 
the two pseudo charts
as shown in Fig.~\ref{fig08}$($a$),($b$)$.
\end{lemma}

\begin{figure}
\centerline{\includegraphics{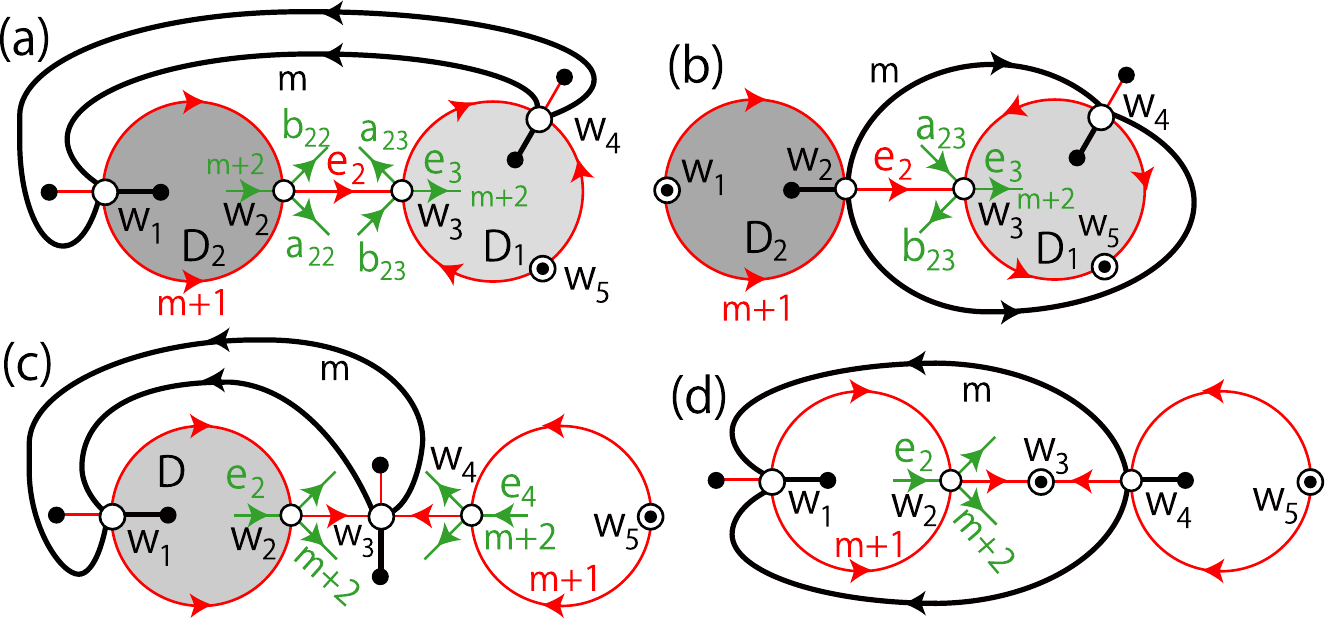}}
\caption{\LABEL{fig08}
(a),(b) Pseudo charts containing the graph as shown in
Fig.~\ref{fig05}(a).
(c),(d) Pseudo charts containing the graph as
shown in Fig.~\ref{fig05}(b). 
}
\end{figure}

\begin{lemma}
{\rm (\cite[Lemma 8.5]{ChartAppVI})}
\LABEL{LemmaTypeH}
Let $\Gamma$ be a minimal chart 
of type $(m;2,3,2)$.
If $\Gamma_{m+1}$ contains the graph as shown in 
Fig.~\ref{fig05}$($b$)$,
then 
$\Gamma$
contains one of the RO-families of the two pseudo charts
as shown in Fig.~\ref{fig08}$($c$),($d$)$.
\end{lemma}

\begin{lemma}
\LABEL{LemmaTypeGTypeHGammaM+2}
Let $\Gamma$ be a minimal chart 
of type $(m;2,3,2)$.
If $\Gamma_{m+2}$ contains one of the two graphs as shown in 
Fig.~\ref{fig05},
then 
$\Gamma$
contains one of the RO-families of the four pseudo charts
as shown in Fig.~\ref{fig09}.
\end{lemma}

\begin{figure}
\centerline{\includegraphics{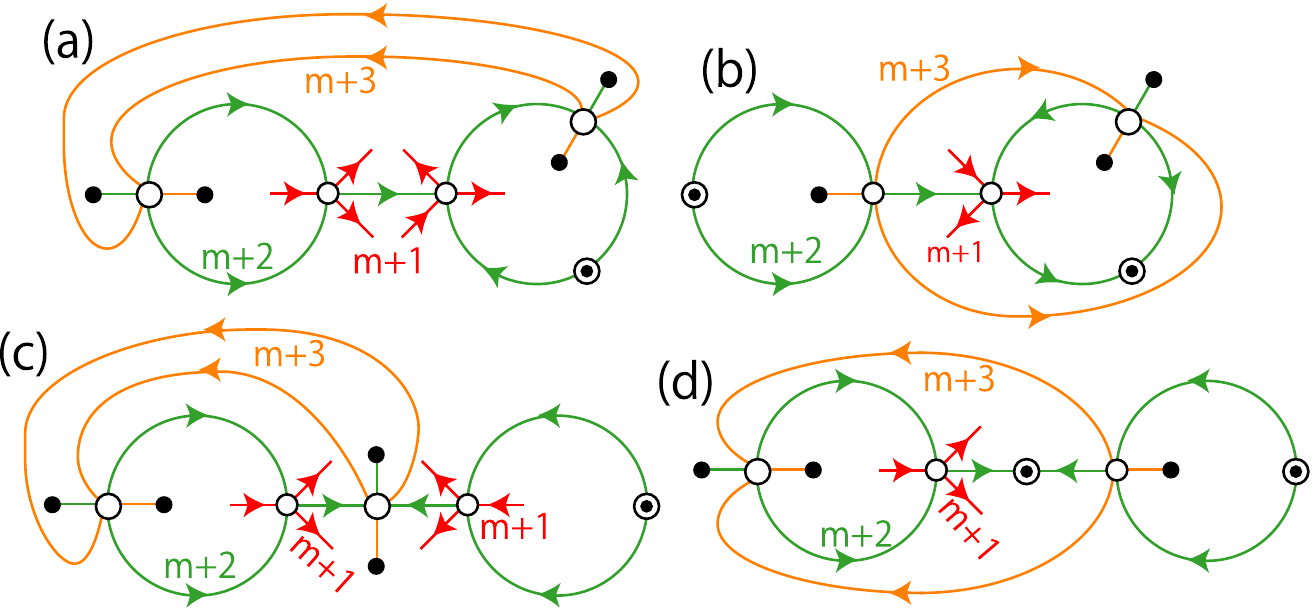}}
\caption{\LABEL{fig09}
(a),(b) Pseudo charts containing the graph as
shown in Fig.~\ref{fig05}(a). 
(c),(d) Pseudo charts containing the graph as
shown in Fig.~\ref{fig05}(b). 
}
\end{figure}

\begin{Proof}
Let $\Gamma$ be a minimal chart of type $(m;2,3,2)$
such that $\Gamma_{m+2}$ contains
the graph as shown in Fig.~\ref{fig05}(a).
Let $\Gamma'$ be the chart obtained from $\Gamma$
by changing labels $\cdots,m,m+1,m+2,m+3,\cdots$
into $\cdots,m+3,m+2,m+1,m,\cdots$,
respectively.
Then $\Gamma'$ is a minimal chart of type $(m;2,3,2)$
such that
 $\Gamma_{m+1}'$ contains
the graph as shown in Fig.~\ref{fig05}(a).
Thus by Lemma~\ref{LemmaTypeG},
the chart $\Gamma'$ contains 
one of the RO-families of the two pseudo charts
as shown in Fig.~\ref{fig08}(a),(b).
Hence $\Gamma$ contains 
one of the RO-families of the two pseudo charts
as shown in Fig.~\ref{fig09}(a),(b).

Similarly we can show this lemma for the case
that $\Gamma_{m+2}$ contains the graph as shown in 
Fig.~\ref{fig05}(b).
\end{Proof}

%



\section{Rings}
\label{s:Ring}

In this section,
we give two examples of non minimal charts $\Gamma$ of
type $(m;2,3,2)$
such that $\Gamma_{m+1}$ contains an oval.

\begin{lemma}
$(${\rm \cite[Lemma 6.1]{ChartAppIII}}$)$
\LABEL{ConditionRing} 
Let $\Gamma$ be a minimal chart.
Let $C$ be a ring or a non simple hoop, 
and $D$ a disk with $\partial D=C$.
If $w(\Gamma\cap D)=1$, then $\Gamma$ is C-move equivalent to the minimal chart $Cl(\Gamma-C)$. 
\end{lemma}

Let $\Gamma$ be a chart, $m$ a label of $\Gamma$, $D$ a disk with $\partial D\subset \Gamma_m$, 
and $k$ a positive integer.
If $\partial D$ contains exactly
$k$ white vertices, 
then $D$ is called 
{\it a $k$-angled disk of $\Gamma_m$}. 
Note that 
the boundary $\partial D$ may contain crossings.

Let $\Gamma$ be a chart, and
$m$ a label of $\Gamma$.
An edge of label $m$ is called a {\it feeler} of a $k$-angled disk $D$ of $\Gamma_m$
if the edge intersects $N-\partial D$
where $N$ is a regular neighborhood of $\partial D$ in $D$.

\begin{lemma}
\LABEL{Theorem2AngledDisk}
{\rm (\cite[Corollary 5.8]{ChartAppII})}
Let $\Gamma$ be a minimal chart.
Let $D$ be a $2$-angled disk of $\Gamma_m$ with at most one feeler.
If $w(\Gamma\cap{\rm Int}D)=0$,
then a regular neighborhood of $D$ contains one of two pseudo charts as shown in Fig.~\ref{fig10}.
\end{lemma}

\begin{figure}
\centerline{\includegraphics{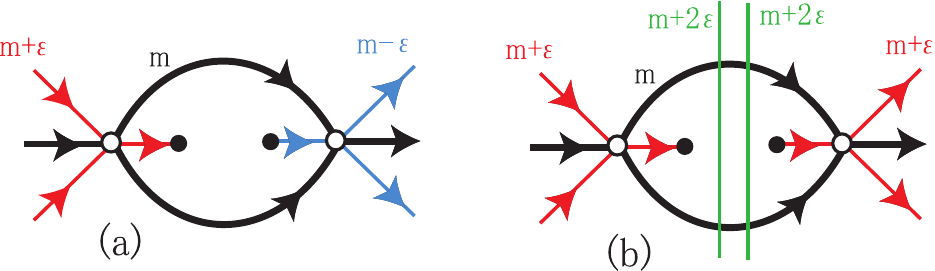}}
\caption{\LABEL{fig10}
$m$ is a label,
and $\varepsilon\in\{+1,-1\}$.}
\end{figure}

Let $\Gamma$ be a chart, 
and $m$ a label of $\Gamma$.
Let $L$ be the closure of a connected component 
of the set obtained by taking out 
all the white vertices from $\Gamma_m$.
If $L$ contains at least one white vertex
but does not contain any black vertex,
then $L$ is called an {\it internal edge of label $m$}.
Note that an internal edge may contain a crossing of $\Gamma$.

A simple arc $\alpha$ in a compact surface $F$ 
is said to be {\it proper} 
if $\alpha\cap\partial F=\partial\alpha$.

\begin{lemma}
\LABEL{NoMinimalOvalGammaM+1}
Let $\Gamma$ be a chart of 
type $(m;2,3,2)$.
If $\Gamma$ contains one of the two pseudo charts
as shown in Fig.~\ref{fig11}$($a$)$,$($b$)$,
then $\Gamma$ is not minimal.
\end{lemma}

\begin{figure}[htb]
\centerline{\includegraphics{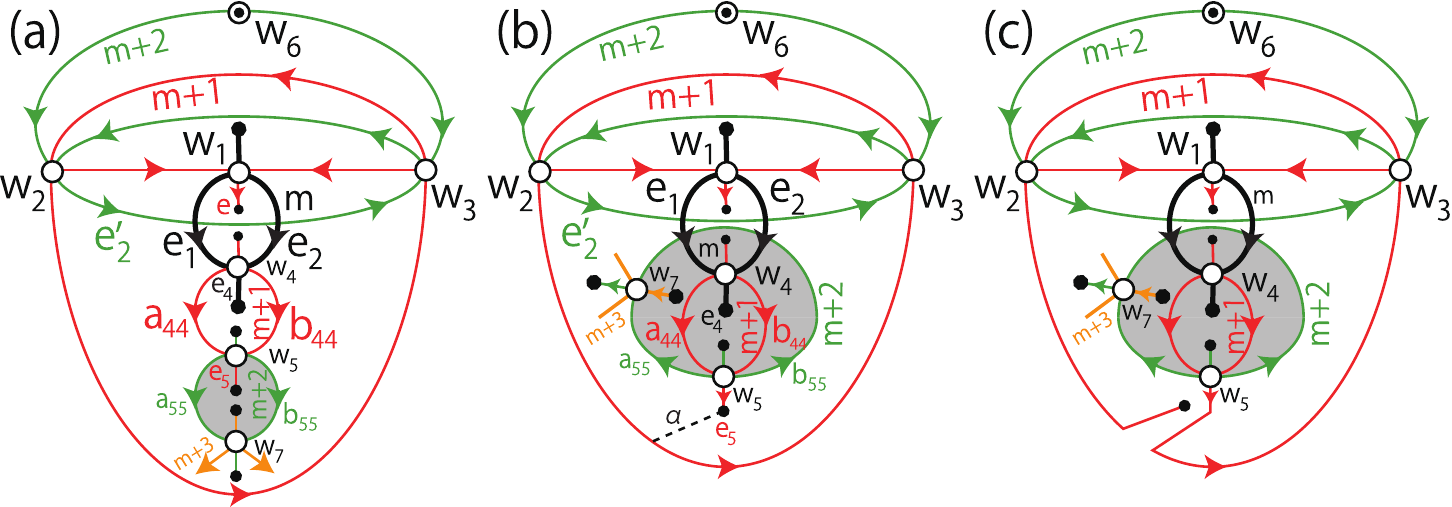}}
\caption{ \LABEL{fig11}
(a) The gray region $F_2$ does not contain $w_4$.
(b) The gray region $F_2$ contains $w_4$.
(c) The terminal edge of label $m+1$ at $w_2$
is not middle at $w_2$.
}
\end{figure}

\begin{Proof}
Suppose that $\Gamma$ is minimal.
We use the notations as shown 
in Fig.~\ref{fig11}(a),(b),
where 
$a_{44},b_{44}$ are internal edges of label $m+1$ at $w_4$,
$e_2'$ is the internal edge of label $m+2$ 
connecting $w_2$ and $w_3$, and
$a_{55},b_{55}$ are internal edges of label $m+2$ at $w_5$.

Since $w_1,w_2,w_3,w_4,w_5\in\Gamma_{m+1}$, 
$w_6,w_7\in\Gamma_{m+2}$ and 
$\Gamma$ is of type $(m;2,3,2)$,
we have
\begin{enumerate}
\item[(1)] $\Gamma_{m+3}$ contains exactly 
two white vertices $w_6,w_7$.
\end{enumerate}

Let $D_1$ be the 3-angled disk of $\Gamma_{m+1}$
with $\partial D_1\ni w_1,w_2,w_3$
and $D_1\ni w_4$.

{\bf Claim~1.} 
There is no ring of 
label $m+2$ in $D_1$.

{\it Proof of Claim~$1$.}
Suppose that there exists a ring $R$ of label 
$m+2$ in $D_1$.
Then the ring $R$ does not intersect 
the two internal edges $a_{44},b_{44}$ of label $m+1$
and the three internal edges $e_2',a_{55},b_{55}$
of label $m+2$.
Thus the ring $R$ is contained in a connected component of
${\rm Int}D_1-(a_{44}\cup b_{44}\cup e_2'\cup a_{55}\cup b_{55})$,
say $F$.
Since each connected component of 
${\rm Int}D_1-(a_{44}\cup b_{44}\cup e_2'\cup a_{55}\cup b_{55})$
does not contain any white vertex,
the ring $R$ does not bound a disk in $F$ 
by Assumption~\ref{Ring}.
Hence $F$ is the open annulus with 
$\partial F\supset e_2'\cup a_{55}\cup b_{55}$.

Since $w_7\in\partial F$ and $w_6\not\in D_1$,
 the two white vertices $w_6,w_7$ 
of $\Gamma_{m+3}$
is separated by the ring $R$ of label $m+2$.
Thus by (1),
there exists a connected component of $\Gamma_{m+3}$
with only one white vertex.
This contradicts Lemma~\ref{LemmaWithTerminal3}.
Hence there is no ring of label $m+2$ in $D_1$.
Hence Claim~1 holds.
{\hfill {$\square$}\vspace{1.5em}}

{\bf Claim~2.} 
By C-moves without changing type of charts, 
we can assume that there is no ring of 
label $m,m+1$ in $D_1$.

{\it Proof of Claim~$2$.}
Suppose that there exists a ring $R$ of label $m$ or $m+1$ 
 in $D_1$.
Then by the similar way of the proof of Claim~1,
the ring $R$ bounds a disk containing 
only one white vertex $w_7$.
Hence by Lemma~\ref{ConditionRing},
the chart $\Gamma$ is C-move equivalent to $Cl(\Gamma-R)$.
Thus we can assume that there is no ring of 
label $m$, $m+1$ in $D_1$.
Hence Claim~2 holds.
{\hfill {$\square$}\vspace{1.5em}}

Now we are ready to prove Lemma~\ref{NoMinimalOvalGammaM+1}.
Suppose that $\Gamma$ contains the pseudo chart
as shown in Fig.~\ref{fig11}(a).
Let $D$ be the 2-angled disk of $\Gamma_m$ in $D_1$
with $\partial D\ni w_1,w_2$.
Then $e_2'\cap D$ is a proper arc of $D$ of label $m+2$.
By Lemma~\ref{Theorem2AngledDisk},
a regular neighborhood of $D$
contains the pseudo chart as shown in Fig.~\ref{fig10}(b).
Thus there exists another proper arc of $D$
contained in a ring of label $m+2$.
This contradicts Claim~1.

Suppose that $\Gamma$ contains the pseudo chart
as shown in Fig.~\ref{fig11}(b).
Let $e_5$ be the terminal edge of label $m+1$ at $w_5$.
Let $\alpha$ be an arc in $D_1$ connecting the black vertex of $e_5$ and a point in $\partial D_1$
(see Fig.~\ref{fig11}(b)).
By Claim~1 and Claim~2,
we can assume that
$(\Gamma_m\cup\Gamma_{m+1}\cap\Gamma_{m+2})\cap
{\rm Int}\alpha=\emptyset$.
Thus we can apply C-II moves and 
a C-I-M2 move between the edge $e_5$ and $\partial D_1$
along the arc $\alpha$
so that we obtain a new terminal edge 
of label $m+1$ at $w_2$
not middle at $w_2$
(see Fig.~\ref{fig11}(c)).
This contradicts Assumption~\ref{AssumeTerminal}.

Therefore $\Gamma$ is not minimal.
We complete the proof of Lemma~\ref{NoMinimalOvalGammaM+1}.
\end{Proof}


\section{IO-Calculation}
\LABEL{s:IOC}

In this section,
we review a useful technic, called IO-Calculation.

Let $\Gamma$ be a chart,
 and $v$ a vertex. 
Let $\alpha$ be a short arc of $\Gamma$ in a small neighborhood of $v$ with $v\in \partial \alpha$. 
If the arc $\alpha$ is oriented to $v$, then $\alpha$ is called {\it an inward arc}, 
and otherwise $\alpha$ is called {\it an outward arc}.

Let $\Gamma$ be an $n$-chart. 
Let $F$ be a closed domain with $\partial F\subset \Gamma_{k-1}\cup\Gamma_{k}\cup \Gamma_{k+1}$ for some label $k$ of $\Gamma$, where $\Gamma_0=\emptyset$ and $\Gamma_{n}=\emptyset$. 
By Condition (iii) for charts,
in a small neighborhood of each white vertex, there are three inward arcs and three outward arcs.
Also in a small neighborhood of each black vertex, there exists only one inward arc or one outward arc.
We often use the following fact, 
when we fix (inward or outward) arcs 
near white vertices and black vertices: 
\begin{enumerate}
\item[$(*)$]
{\it The number of inward arcs contained in $F\cap \Gamma_k$ is equal to the number of outward arcs in $F\cap \Gamma_k$.
}
\end{enumerate}
When we use this fact, 
we say that we use {\it IO-Calculation with respect to $\Gamma_k$ in $F$}.
For example, in a minimal chart $\Gamma$, 
consider the pseudo chart as shown in Fig.~\ref{fig12} 
where
\begin{enumerate}
\item[(1)] $F$ is an annulus with $\partial F\subset\Gamma_{k-1}$,
\item[(2)]  $v_1,v_2,v_3,v_4,v_5$ are white vertices 
in $\partial F$  with  
$v_1,v_4\in\Gamma_{k-2}\cap\Gamma_{k-1}$ and
$v_2,v_3,v_5\in\Gamma_{k-1}\cap\Gamma_{k}$,
\item[(3)] $e_2$ is an edge of label $k$ 
oriented outward at $v_2$,
\item[(4)] $e_3$ is an edge of label $k$ 
oriented inward at $v_3$ but not middle at $v_3$,
\item[(5)] $a_{55},b_{55}$ are edges of label $k$ 
oriented outward at $v_5$ 
but not middle at $v_5$.
\end{enumerate}
Then we can show that $w(\Gamma\cap{\rm Int}F)\ge1$.
Suppose $w(\Gamma\cap{\rm Int}F)=0$.
By (4) and (5), 
none of $e_3,a_{55},b_{55}$ are middle at $v_3$ or $v_5$.
Thus by Assumption~\ref{AssumeTerminal},
\begin{enumerate}
\item[(6)] none of $e_3,a_{55},b_{55}$ are terminal edges. 
\end{enumerate}
If $e_2$ is a terminal edge, then
the number of inward arcs in $F\cap \Gamma_{k}$ is two,  
but the number of outward arcs in $F\cap \Gamma_{k}$ is three. 
This contradicts the fact $(*)$.
If $e_2$ is not a terminal edge,
then we have the same contradiction. 
Thus $w(\Gamma\cap{\rm Int}F)\ge1$.
Instead of the above argument, 
we just say that 
\begin{enumerate}
\item[]
{\it we have $w(\Gamma\cap{\rm Int}F)\ge1$ 
by IO-Calculation with respect to $\Gamma_{k}$ in $F$.}
\end{enumerate}

\begin{figure}[htb]
\centerline{\includegraphics{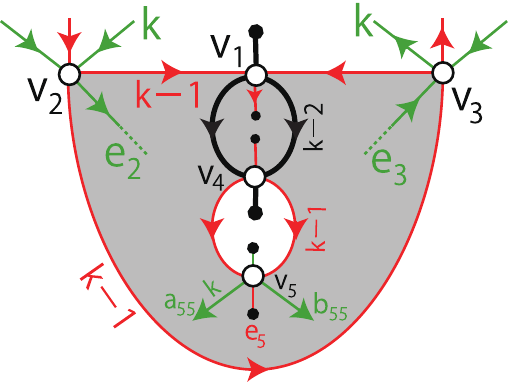}}
\caption{\LABEL{fig12} The gray region is the annulus $F$.}
\end{figure}

\begin{lemma}$($Triangle Lemma$)$
{\rm (\cite[Lemma 8.3(2)]{ChartAppIV})}
\LABEL{LemmaTriangle}
 For a minimal chart $\Gamma$, 
if there exists a $3$-angled disk $D_1$ of $\Gamma_m$ without feelers in a disk $D$ as shown in Fig.~\ref{fig13}$($a$)$,
then $w(\Gamma\cap${\rm Int}$D_1)\ge1$.
\end{lemma}

\begin{figure}
\centerline{\includegraphics{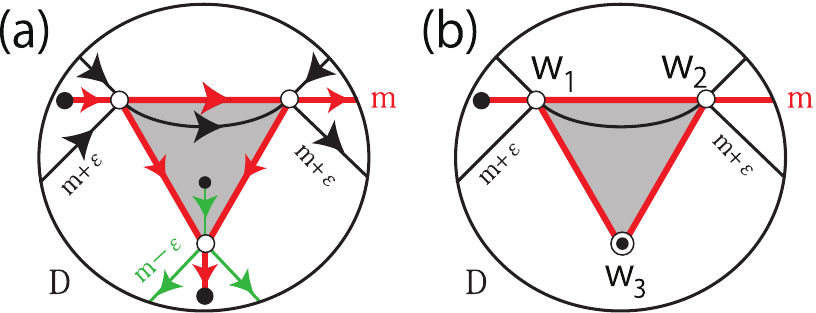}}
\caption{\LABEL{fig13}
The gray region is the 3-angled disk $D_1$. 
The thick lines are edges of label $m$,
and $\varepsilon\in\{+1,-1\}$.}
\end{figure}

\begin{corollary}$($Corollary of Triangle Lemma$)$
\LABEL{CorollaryTriangleLemma}
Let $\Gamma$ be a minimal chart. 
Let  $D_1$ be a $3$-angled disk of $\Gamma_m$ 
with at most one feeler in a disk $D$ as shown 
in Fig.~\ref{fig13}$($b$)$.
Let $w_1,w_2,w_3$ be the white vertices in $\partial D_1$
with $w_1,w_2\in\Gamma_{m+\varepsilon}$ 
$(\varepsilon\in\{+1,-1\})$.
If $w_3\in\Gamma_{m-\varepsilon}$,
then $w(\Gamma\cap${\rm Int}$D_1)\ge1$.
\end{corollary}

\begin{Proof}
Let $e$ be the terminal edge of label $m$ at $w_3$.
Without loss of generality we can assume
that the edge $e$ is oriented outward at $w_3$.

If $e\not\subset D_1$,
then the 3-angled disk $D_1$ is as shown in
Fig.~\ref{fig13}(a).
Thus by Lemma~\ref{LemmaTriangle},
we have $w(\Gamma\cap{\rm Int}D_1)\ge1$.

If $e\subset D_1$,
then there exist two internal edges $e',e''$
(possibly terminal edges)
of label $m-\varepsilon$ at $w_3$ in $D_1$.
By Assumption~\ref{AssumeTerminal},
neither $e'$ nor $e''$
is a terminal edge,
and both of $e',e''$ are oriented outward at $w_3$.
Thus we have $w(\Gamma\cap{\rm Int}D_1)\ge1$
by IO-Calculation with respect to 
$\Gamma_{m-\varepsilon}$ in $D_1$.
\end{Proof}


\section{X-change Lemma}
\label{s:XchangeLemma}

In this section
we investigate a minimal chart $\Gamma$
of type $(m;2,3,2)$
such that $\Gamma_{m+1}$ contains an oval.
We shall show that 
the chart $\Gamma$ contains the pseudo chart
as shown in Fig.~\ref{fig17}(a)
by C-moves without changing types of charts.

Let $\Gamma$ and $\Gamma^\prime $ be C-move equivalent charts. 
Suppose that a pseudo chart $X$ of $\Gamma$ is also a pseudo chart of $\Gamma^\prime$. 
Then we say that 
$\Gamma$ is modified to $\Gamma^\prime$ by {\it C-moves keeping $X$ fixed}.
In Fig.~\ref{fig14},
we give examples of C-moves keeping pseudo charts  fixed.

\begin{figure}[htb]
\centerline{\includegraphics{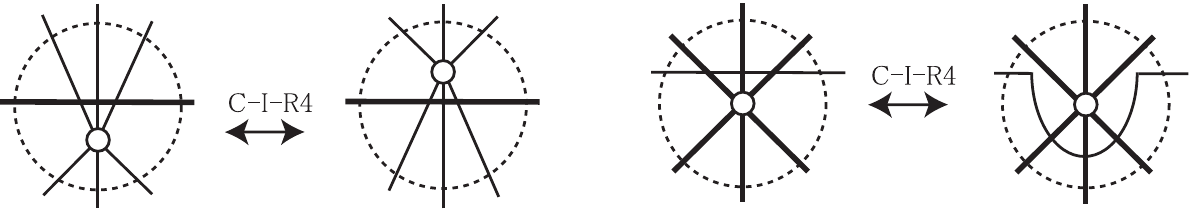}}
\caption{\LABEL{fig14} 
C-moves keeping thicken figures fixed.}
\end{figure}

Let $\Gamma$ be a chart and $m$ a label.
An oval $G$ of $\Gamma_{m+1}$
 is said to be {\it special},
if there exists a 2-angled disk $D$ of $\Gamma_{m+1}$ without feelers
such that $\partial D\subset G$, 
$w(\Gamma\cap{\rm Int}D)=0$, 
the disk $D$ contains a terminal edge of label $m$ and a terminal edge of label $m+2$, but  $D$ does not contain any free edges, hoops nor crossings  
(see Fig.~\ref{fig15}(a)).


\begin{figure}
\centerline{\includegraphics{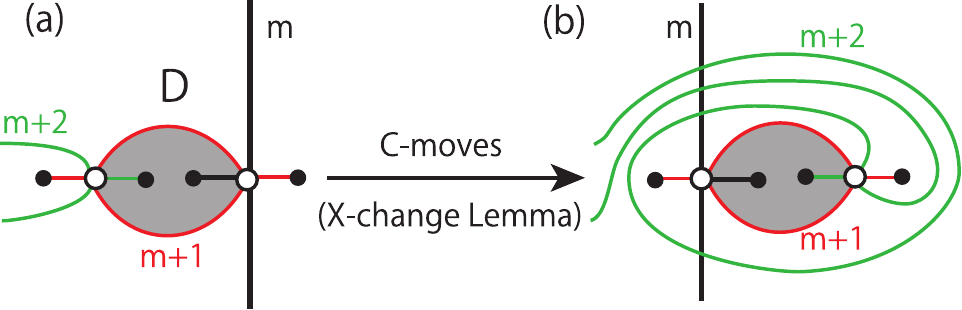}}
\caption{\LABEL{fig15}
The gray regions are the disk $D$.}
\end{figure}

\begin{lemma}
{\rm (\cite[Lemma 6.1 and Lemma 6.3]{INS})}
\LABEL{LemmaTwist}
Let $\Gamma$ be a chart.
Let $G$ be an oval of $\Gamma_{m+1}$
and $D$ a $2$-angled disk of $\Gamma_{m+1}$ without feelers such that $\partial D\subset G$ and 
$w(\Gamma\cap{\rm Int}D)=0$.
\begin{enumerate}
\item[{\rm (a)}] $($X-change Lemma$)$
If $G$ is a special oval in 
a minimal chart $\Gamma$,
then the chart $\Gamma$ is C-move equivalent to the chart obtained from $\Gamma$ by replacing a regular neighborhood of $D$ with the pseudo chart as shown in Fig.~\ref{fig15}$($b$)$.
\item[{\rm (b)}]
 If $D$ contains a terminal edge of label $m$ and a terminal edge of label $m+2$,
then $G$ can be modified to a special oval by C-moves in a regular neighborhood of $D$ keeping $G\cup\Gamma_m\cup\Gamma_{m+2}$ fixed.
\end{enumerate}
\end{lemma}

\begin{lemma}
\LABEL{OvalOneFeeler}
{\rm (\cite[Lemma 4.4(2)]{ChartAppV})}
Let $\Gamma$ be a minimal chart,
and $m$ a label of $\Gamma$.
Let $D$ be a $2$-angled disk of $\Gamma_m$
with exactly one feeler. 
 If $w(\Gamma\cap {\rm Int}D)=1$,
then a regular neighborhood of $D$ contains 
one of the RO-family of
the pseudo chart as shown in 
Fig.~\ref{fig16}.
\end{lemma}

\begin{figure}[htb]
\centerline{\includegraphics{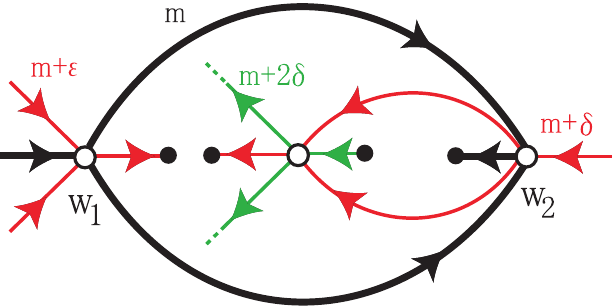}}
\caption{\LABEL{fig16}
$m$ is a label, $\varepsilon,\delta\in\{+1,-1\}$.}
\end{figure}


\begin{lemma}
\LABEL{OvalGammaM+1Step1A}
Let $\Gamma$ be a minimal chart 
of type $(m;2,3,2)$.
If $\Gamma$ contains the pseudo chart
as shown in Fig.~\ref{fig07}$($a$)$,
then the chart $\Gamma$ is C-move equivalent to
 a minimal chart of type $(m;2,3,2)$
containing the pseudo chart
as shown in Fig.~\ref{fig17}$($a$)$.
Moreover
$e_2'=e_3'$ and
$a_{55}\cap b_{55}$ contains a white vertex
different from $w_5$,
where $e_2',e_3',a_{55},b_{55}$ are internal edges
of label $m+2$
at $w_2,w_3,w_5,w_5$ respectively.
\end{lemma}

\begin{figure}
\centerline{\includegraphics{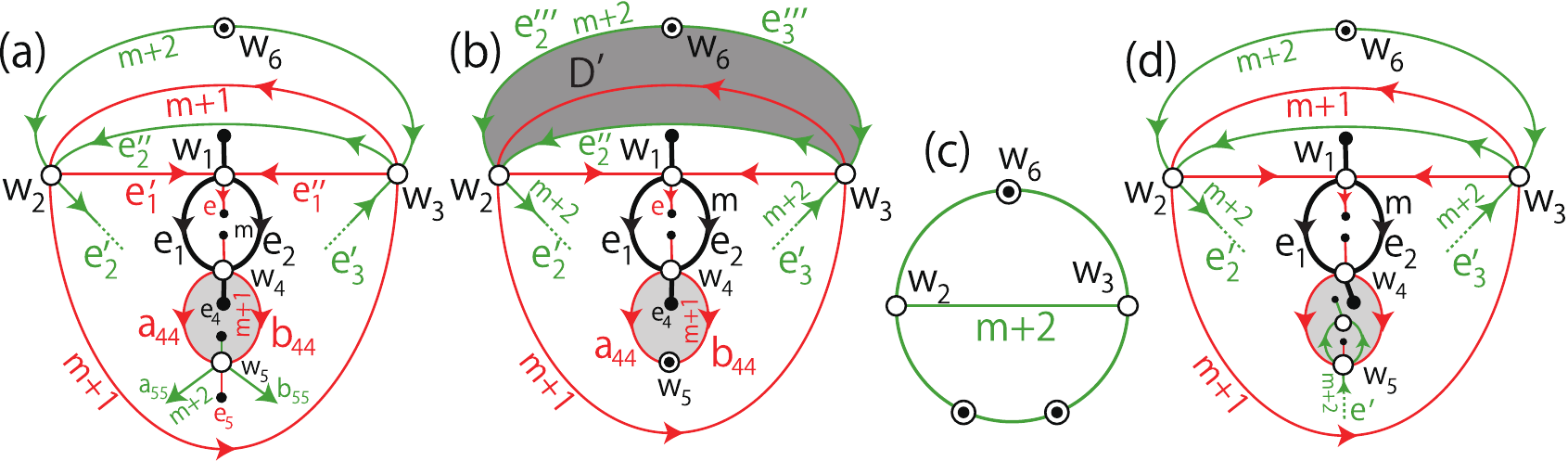}}
\caption{
 The light gray regions are the disk $E$.
 The dark gray region is the disk $D'$.
 \LABEL{fig17}
}
\end{figure}

\begin{Proof}
We use the notations as shown in 
Fig.~\ref{fig07}(a).
Now $\Gamma_{m+1}$ contains a skew $\theta$-curve
(see Fig.~\ref{fig07}(a)).
The skew $\theta$-curve separates $S^2$
into three disks.
Let $D_1,D_2,D_3$ be the three disks
such that
$D_1\supset a_{44}\cup b_{44}$, 
the disk $D_2$ contains the terminal edge 
of label $m$ at $w_1$,
and
$D_3$ is a 2-angled disk.

Let $e_2',e_2'',e_2'''$ be internal edges 
(possibly terminal edges) of label $m+2$ at $w_2$
in $D_1,D_2,D_3$ respectively.
Let $e_3',e_3'',e_3'''$ be internal edges 
(possibly terminal edges) of label $m+2$ at $w_3$
in $D_1,D_2,D_3$ respectively 
(see Fig.~\ref{fig07}(a)).
Then
\begin{enumerate}
\item[(1)] the two edges $e'''_{2},e'''_{3}$
are oriented inward at $w_2,w_3$
respectively,
\item[(2)] the edge $e_2'$ is oriented outward at $w_2$.
\end{enumerate}

{\bf Claim 1.}
$w(\Gamma\cap{\rm Int}D_3)\ge1$.

{\it Proof of Claim~$1$.}
We can show that
neither $e_2'''$ nor $e_3'''$
is middle at $w_2$ or $w_3$.
Thus by Assumption~\ref{AssumeTerminal},
\begin{enumerate}
\item[(3)] neither $e'''_{2}$ nor $e'''_{3}$
is a terminal edge.
\end{enumerate}
Hence (1) implies that
$w(\Gamma\cap{\rm Int}D_3)\ge1$
by IO-Calculation with respect to $\Gamma_{m+2}$
in $D_3$.
Thus Claim~1 holds.
{\hfill {$\square$}\vspace{1.5em}}

We shall show $w_5\in\Gamma_{m+1}\cap\Gamma_{m+2}$.
Since $w_1,w_2\in\Gamma_m\cap\Gamma_{m+1}$,
$w_5\in\Gamma_{m+1}$ and since
$\Gamma$ is of type $(m;2,3,2)$,
we have $w_5\in\Gamma_{m+1}\cap\Gamma_{m+2}$.

Let $e_5$ be the terminal edge of label $m+1$ at $w_5$, 
and
$a_{55},b_{55}$ the internal edges (possibly terminal edges) of label $m+2$ at $w_5$
such that $a_{55},e_{5},b_{55}$ are consecutive edges 
at $w_5$.
Since  $a_{44},b_{44}$ are oriented from $w_4$ to $w_5$,
by Assumption~\ref{AssumeTerminal}
\begin{enumerate}
\item[(4)] the two edges $a_{55},b_{55}$
are oriented outward at $w_5$,
\item[(5)] neither $a_{55}$ nor $b_{55}$
is a terminal edge.
\end{enumerate}

Let $E$ be the 2-angled disk of $\Gamma_{m+1}$ in $D_1$
with $\partial E=a_{44}\cup b_{44}$.

{\bf Claim 2.}
If $e_5\subset E$, 
then $w(\Gamma\cap{\rm Int}E)\ge1$,
otherwise $w(\Gamma\cap({\rm Int}D_1-E))\ge1$.

{\it Proof of Claim~$2$.}
If $e_5\subset E$,
then (4) and (5) imply that 
$w(\Gamma\cap{\rm Int}E)\ge1$
by IO-Calculation with respect to $\Gamma_{m+2}$ in $E$.

If $e_5\not\subset E$, then
considering as $F=Cl(D_1-E)$ and $k=m+2$ 
in the example of IO-Calculation in Section~\ref{s:IOC},
we have $w(\Gamma\cap({\rm Int}D_1-E))\ge1$.
Hence Claim~2 holds.
{\hfill {$\square$}\vspace{1.5em}}

{\bf Claim 3.}
$w(\Gamma\cap{\rm Int}D_1)=3$,
$w(\Gamma\cap{\rm Int}D_2)=0$, and
$w(\Gamma\cap{\rm Int}D_3)=1$.

{\it Proof of Claim~$3$.}
Since $E\subset {\rm Int}D_1$,
we have $w(\Gamma\cap{\rm Int}D_1)\ge3$
by Claim~2.
By Claim~1,
we have $w(\Gamma\cap{\rm Int}D_3)\ge1$.

Suppose that
$w(\Gamma\cap{\rm Int}D_1)>3$ or
$w(\Gamma\cap{\rm Int}D_2)>0$ or
$w(\Gamma\cap{\rm Int}D_3)>1$.
Since $\Gamma$ is of type $(2,3,2)$,
we have $w(\Gamma)=7$.
Since the skew $\theta$-curve in $\Gamma_{m+1}$
contains exactly three white vertices,\\
$7=w(\Gamma)=3+w(\Gamma\cap{\rm Int}D_1)+
w(\Gamma\cap{\rm Int}D_2)+
w(\Gamma\cap{\rm Int}D_3)>3+3+0+1=7$.\\
This is a contradiction.
Hence Claim~3 holds.
{\hfill {$\square$}\vspace{1.5em}}

{\bf Claim 4.}
$e_2''=e_3''$ and 
$e_2'''\cap e_3'''$ contains a white vertex.
Moreover 
the chart $\Gamma$ contains the pseudo chart
as shown in Fig.~\ref{fig17}(b).

{\it Proof of Claim~$4$.}
We can show that 
the edge $e_2''$ is not a terminal edge
by Assumption~\ref{AssumeTerminal}.
Since $w(\Gamma\cap{\rm Int}D_2)=0$ by Claim~3,
we have $e_2''=e_3''$.

By (1) and (3),
both edges $e_2''',e_3'''$ contain white vertices $w',w''$ 
different from
$w_2,w_3$
respectively.
Since $w(\Gamma\cap{\rm Int}D_3)=1$ by Claim~3,
the white vertices $w',w''$ are the same vertex.
Hence $e_2'''\cap e_3'''$ contains a white vertex, 
say $w_6$.
Moreover there exists a terminal edge 
of label $m+2$ at $w_6$,
i.e. $w_6$ is a BW-vertex with respect to $\Gamma_{m+2}$.
Hence
the chart $\Gamma$ contains the pseudo chart
as shown in Fig.~\ref{fig17}(b).
Thus Claim~4 holds.
{\hfill {$\square$}\vspace{1.5em}}

{\bf Claim~5.}
The edge $e_2'$ is not a terminal edge.

{\it Proof of Claim~$5$.}
Let $D'$ be the 3-angled disk of $\Gamma_{m+2}$
with $\partial D'=e_2''\cup e_2'''\cup e_3'''$
and $w(\Gamma\cap {\rm Int}D')=0$
(see Fig.~\ref{fig17}(b)).

If $e_2'$ is a terminal edge,
then $w(\Gamma\cap {\rm Int }D')\ge1$
by Triangle Lemma (Lemma~\ref{LemmaTriangle}).
This is a contradiction.
Hence edge $e_2'$ is not a terminal edge.
Thus Claim~5 holds.
{\hfill {$\square$}\vspace{1.5em}}

By Claim~2 and the first equation of Claim~3,
we have the following claim:

{\bf Claim 6.}
If $e_5\subset E$, 
then $w(\Gamma\cap{\rm Int}E)=1$ and 
$w(\Gamma\cap({\rm Int}D_1-E))=0$,
otherwise $w(\Gamma\cap({\rm Int}D_1-E))=1$.

If $e_5\subset E$,
then let $w_7$ be the white vertex in ${\rm Int}E$,
otherwise 
let $w_7$ be the white vertex in ${\rm Int}D_1-E$.

Now we are ready to prove Lemma~\ref{OvalGammaM+1Step1A}.
There are two cases:
(i) $e_5\not\subset E$,
(ii) $e_5\subset E$.

{\bf Case (i).} 
By $e_5\not\subset E$,
the chart $\Gamma$ contains the pseudo chart
as shown in Fig.~\ref{fig17}(a).
We shall show that
$e_2'=e_3'$ and $a_{55}\cap b_{55}\ni w_7$.

By Claim~5, the edge $e_2'$ is not a terminal edge.
Thus  by (2),(4) and Claim~6,
we have $e_2'=e_3'$ or $e_2'\ni w_7$. 
Similarly
by (2),(4),(5) and Claim~6,
we can show that
($a_{55}=e_3'$ or $a_{55}\ni w_7$) and
($b_{55}=e_3'$ or $b_{55}\ni w_7$).
Therefore
($e_2'=e_3'$, $a_{55}\ni w_7$, $b_{55}\ni w_7$) or
($e_2'\ni w_7$, $a_{55}=e_3'$, $b_{55}\ni w_7$) or
($e_2'\ni w_7$, $a_{55}\ni w_7$, $b_{55}=e_3'$).

For the second case and third case,
the graph $\Gamma_{m+2}$ contains the graph as shown in
Fig.~\ref{fig17}(c).
However by Lemma~\ref{OvalTypeGTypeH},
the graph $\Gamma_{m+2}$ does not contain
the graph as shown in Fig.~\ref{fig17}(c).
Thus the second case and third case do not occur.
Hence $e_2'=e_3'$, $a_{55}\cap b_{55}\ni w_7$.

{\bf Case (ii).}
By Claim~6, we have $w(\Gamma\cap{\rm Int}E)=1$.
Thus by Lemma~\ref{OvalOneFeeler},
the disk $E$ contains the pseudo chart as shown in
Fig.~\ref{fig16} 
(see Fig.~\ref{fig17}(d)).
Hence $a_{55}\cap b_{55}\ni w_7$.

Now, the edge $e_3'$ is oriented inward at $w_3$
but not middle at $w_3$.
Thus by Assumption~\ref{AssumeTerminal},
the edge $e_3'$ is not a terminal edge.
Hence by Claim~5,
we have $e_2'=e_3'$.

Let $e'$ be an internal edge (possibly terminal edge)
of label $m+2$ at $w_5$ in $Cl(D_1-E)$.
Since $w(\Gamma\cap({\rm Int}D_1-E))=0$ by Claim~6,
the condition $e_2'=e_3'$ implies that
the edge $e'$ is a terminal edge.

Hence by  Lemma~\ref{LemmaTwist}(b)
there exists a special oval of label $m+2$
by C-moves without changing the type of charts.
Thus by X-change Lemma (Lemma~\ref{LemmaTwist}(a)),
the chart $\Gamma$
is C-move equivalent to a minimal chart of type $(m;2,3,2)$
containing the pseudo chart
as shown in Fig.~\ref{fig17}(a).

Therefore, 
we complete the proof of Lemma~\ref{OvalGammaM+1Step1A}.
\end{Proof}


Similarly by applying several times of 
Lemma~\ref{LemmaTwist},
we can show the following lemma.

\begin{lemma}
\LABEL{OvalGammaM+1Step1B}
Let $\Gamma$ be a minimal chart 
of type $(m;2,3,2)$.
If $\Gamma$ contains the pseudo chart
as shown in Fig.~\ref{fig07}$($b$)$,
then the chart $\Gamma$ is C-move equivalent
to a minimal chart of type $(m;2,3,2)$
containing the pseudo chart
as shown in Fig.~\ref{fig17}$($a$)$.
Moreover
$e_2'=e_3'$ and
$a_{55}\cap b_{55}$ contains a white vertex
different from $w_5$.
\end{lemma}

\section{Disk Lemma}
\LABEL{s:DiskLemma}

In this section
we shall show that
neither $\Gamma_{m+1}$ nor $\Gamma_{m+2}$
contains an oval 
for any minimal chart $\Gamma$ of type $(m;2,3,2)$.

Let $\Gamma$ be a chart, and $D$ a disk.
Let $\alpha$ be a simple arc in $\partial D$,
and $\gamma$ a simple arc in an internal edge of label $k$.
The simple arc $\gamma$ is called
a {\it {$(D,\alpha)$-arc}} of label $k$
provided that 
$\partial \gamma \subset $Int$\alpha$
and
Int$\gamma\subset $Int$D$. 
If there is no $(D,\alpha)$-arc in $\Gamma$,
then the chart $\Gamma$ is said to be
$(D,\alpha)$-{\it arc free}.

\begin{lemma}
$($New Disk Lemma$)$
{\em (\cite[Lemma 7.1]{Chart33},
cf. \cite[Lemma 3.2]{ChartApp1})} 
\LABEL{NewDiskLemma}
Let $\Gamma$ be a chart and
$D$ a disk 
whose interior does not contain 
a white vertex nor a black vertex of $\Gamma$.
Let $\alpha$ be a simple arc in $\partial D$ 
such that ${\rm Int}\alpha$ does not contain 
a white vertex nor a black vertex of $\Gamma$.
Let $V$ be a regular neighborhood of $\alpha$. 
Suppose that 
the arc $\alpha$ is contained in 
an internal edge of some label $k$ of $\Gamma$.
Then by applying C-I-M2 moves, C-I-R2 moves, 
and C-I-R3 moves in $V$, 
there exists 
a $(D,\alpha)$-arc free chart $\Gamma'$ 
obtained from the chart $\Gamma$ 
keeping $\alpha$ fixed 
$($see Fig.~\ref{fig18}$)$.
\end{lemma}

\begin{figure}
\centerline{\includegraphics{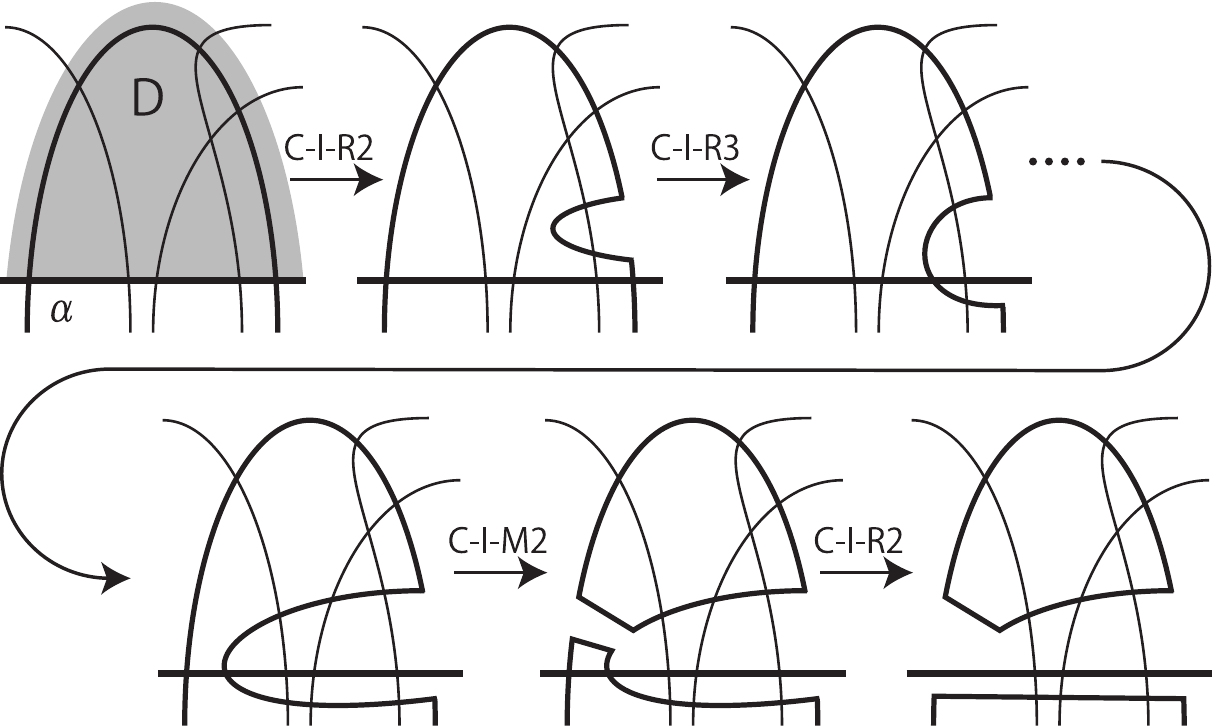}}
\caption{\LABEL{fig18}
The gray region is the disk $D$.}
\end{figure}

\begin{lemma}
\LABEL{EdgeMM-1M+1}
{\rm (\cite[Lemma 8.2]{ChartAppIV})}
Let $\Gamma$ be a chart,
$e$ an internal edge of label $m$,
and $w_1,w_2$ the white vertices in $e$.
Suppose $w_1\in\Gamma_{m-1}$ and $w_2\in\Gamma_{m+1}$.
Then for any neighborhood $V$ of the edge $e$,
there exists a chart $\Gamma'$ obtained 
from the chart $\Gamma$ by C-I-R2 moves, C-I-R3 moves
and C-I-R4 moves in $V$ keeping 
$\Gamma_{m-1}\cup\Gamma_m\cup\Gamma_{m+1}$ fixed
such that the edge $e$ does not contain any crossings.
\end{lemma}


\begin{lemma}
{\rm (\cite[Lemma 7.1]{ChartAppVI})}
\LABEL{GammaM+1ToGammaM+2Cite11}
Let $G$ be one of 12 graphs as shown in
Fig.~3 and Fig.~4 in \cite{ChartAppVI}.
If 
for any minimal chart $\Gamma$ of type $(m;2,3,2)$,
the graph $\Gamma_{m+1}$ does not contain
the graph $G$,
then the graph $\Gamma_{m+2}$ does not contain
the graph $G$.
\end{lemma}

Since the four graphs as shown in
Fig.~\ref{fig04} and Fig.~\ref{fig05}
are four graphs in the 12 graphs 
of the above lemma,
we have the following lemma:

\begin{lemma}
\LABEL{GammaM+1ToGammaM+2}
Let $G$ be one of four graphs as shown in
Fig.~\ref{fig04} and Fig.~\ref{fig05}.
If 
for any minimal chart $\Gamma$ of type $(m;2,3,2)$,
the graph $\Gamma_{m+1}$ does not contain
the graph $G$,
then the graph $\Gamma_{m+2}$ does not contain
the graph $G$.
\end{lemma}

\begin{proposition}
\LABEL{PropOvalGammaM+1}
{\rm (cf. \cite[Lemma 6.3]{ChartAppV})}
Let $\Gamma$ be a minimal chart 
of type $(m;2,3,2)$.
Then neither $\Gamma_{m+1}$ nor $\Gamma_{m+2}$
contains an oval.
\end{proposition}

\begin{Proof}
Suppose that $\Gamma_{m+1}$ contains an oval.
By Lemma~\ref{OvalGammaM+1Step0}, 
Lemma~\ref{OvalGammaM+1Step1A} and 
Lemma~\ref{OvalGammaM+1Step1B},
the chart $\Gamma$ contains the pseudo chart
as shown in Fig.~\ref{fig17}$($a$)$,
where $e_2',e_3',a_{55},b_{55}$ are internal edges 
of label $m+2$ at $w_2,w_3,w_5,w_5$ respectively
and
\begin{enumerate}
\item[(1)] $e_2'=e_3'$,
\item[(2)] $a_{55}\cap b_{55}$ contains a white vertex
different from $w_5$, say $w_7$.
\end{enumerate}

Since $w_1,w_2,w_3,w_4,w_5\in \Gamma_{m+1}$, 
$w_6,w_7\in\Gamma_{m+2}$
and since $\Gamma$ is of type $(m;2,3,2)$,
we have
\begin{enumerate}
\item[(3)] $\Gamma_{m+3}$ contains exactly two white vertices $w_6,w_7$.
\end{enumerate}

By (1),
the union $e_2'\cup e_2''$ bounds a 2-angled disk
of $\Gamma_{m+2}$.
Let $E'$ be the 2-angled disk of $\Gamma_{m+2}$
with
$\partial E'=e_2'\cup e_2''$ and 
${\rm Int}E'\ni w_1$.

{\bf Claim 1.}
${\rm Int}E'$ contains exactly one white vertex $w_1$.

{\it Proof of Claim~$1$.}
It is possible that 
${\rm Int}E'$ contains $w_1,w_4,w_5,w_7$.

We shall show $w_7\not\in {\rm Int}E'$.
If $w_7\in {\rm Int}E'$,
then 
the two white vertices $w_6$ and $w_7$ 
of $\Gamma_{m+3}$
are separated by
the simple closed curve $\partial E'$ 
of label $m+2$.
Thus by (3)
there exists a connected component of $\Gamma_{m+3}$
with exactly one white vertex.
This contradicts Lemma~\ref{LemmaWithTerminal3}.
Hence $w_7\not\in {\rm Int}E'$.

We shall show that neither $w_4$ nor $w_5$ is 
contained in ${\rm Int}E'$.
If one of $w_4,w_5$ is contained in ${\rm Int}E'$,
then
both of $w_4$ and $w_5$ are
contained in ${\rm Int}E'$, 
because 
 $\partial E'\subset \Gamma_{m+2}$ and
 $w_4,w_5\in a_{44}\subset \Gamma_{m+1}$.
Moreover $w_5\in a_{55}\cap b_{55}$ implies
$a_{55}\cup b_{55}\subset {\rm Int}E'$.
Thus $w_7\in a_{55}\cap b_{55}$ by (2) implies
$w_7\in {\rm Int}E'$.
This is a contradiction.
Hence neither $w_4$ nor $w_5$ is
contained in ${\rm Int}E'$.
Therefore ${\rm Int}E'$ contains 
exactly one white vertex $w_1$.
Thus Claim~1 holds.
{\hfill {$\square$}\vspace{1.5em}}

Let $e_1',e_1''$ be the internal edges of label $m+1$
at $w_1$ (see Fig.~\ref{fig17}(a)).
Then the arc $e_1'\cup e_1''$ separates the disk $E'$
into two disks.
Let $F_1$ be one of the two disks 
with $\partial F_1\supset e_2'$.
By Claim~1,
we have
\begin{enumerate}
\item[(4)] $w(\Gamma\cap{\rm Int}F_1)=0$.
\end{enumerate}

Let $e_1,e_2$ be the internal edges of label $m$
connecting $w_1$ and $w_4$.

{\bf Claim 2.}
We can assume that
for each $i=1,2$ 
the intersection $e_i \cap e_2'$ is one point by C-moves
in a neighborhood of $F_1$
keeping $\Gamma_{m+1}\cup \Gamma_{m+2}\cup \Gamma_{m+3}$
fixed (see Fig.~\ref{fig19}(a)).

{\it Proof of Claim~$2$.}
Since $w_4$ is contained in the outside of $E'$
by Claim~1,
for each $i=1,2$
the intersection $e_i\cap E'$ consists of 
an arc with $w_1$ and 
proper arcs of $E'$.
Since $F_1$ is obtained from $E'$ by cutting along
the arc $e_1'\cup e_1''$ in $\Gamma_{m+1}$,
\begin{enumerate}
\item[(5)] $e_i\cap F_1=e_i\cap E'$ consists of an arc with $w_1$
and $(F_1,e_2')$-arcs of label $m$
(see Fig~\ref{fig19}(b)).
\end{enumerate}

Let $N$ be a regular neighborhood of 
the terminal edge $e$ of label $m+1$ at $w_1$ in $F_1$.
Then $\widetilde{F_1}=Cl(F_1-N)$ is a disk.
By (4), we can apply New Disk Lemma 
(Lemma~\ref{NewDiskLemma}) 
for the disk $\widetilde{F_1}$
so that 
the chart $\Gamma$ is $(\widetilde{F_1},e_2')$-arc free.
Thus by (5),
$e_i\cap F_1$ is an arc with $w_1$.
Hence $e_i\cap e_2'$ is one point.
Thus Claim~2 holds.
{\hfill {$\square$}\vspace{1.5em}}

By (3),
we have $w_7\in\Gamma_{m+3}$.
Hence by (2), the edge $a_{55}$ of label $m+2$
connects $w_5\in\Gamma_{m+1}$ and $w_7\in\Gamma_{m+3}$.
Thus by Lemma~\ref{EdgeMM-1M+1}
we can assume that
\begin{enumerate}
\item[(6)] the edge $a_{55}$ does not contain any crossing.
\end{enumerate}

Let $F_2$ be the 2-angled disk of $\Gamma_{m+2}$
with $\partial F_2=a_{55}\cup b_{55}$ and $w_6\not\in F_2$.
Since $w_7\in\partial F_2$,
we have 
\begin{enumerate}
\item[(7)] $w(\Gamma_{m+3}\cap{\rm Int}F_2)=0$ by (3),
\item[(8)] ${\rm Int}F_2$ contains at most 
one white vertex $w_4$.
\end{enumerate}

Let $e_7$ be the terminal edge of label $m+2$ at $w_7$.

{\bf Claim 3.} $e_7\not\subset F_2$.

{\it Proof of Claim~$3$.}
Suppose $e_7\subset F_2$.
Then by (3),
there exist two internal edges $e_7',e_7''$
(possibly terminal edges) of label $m+3$ 
at $w_7$ in $F_2$.
Since both of the edges $a_{55}$ and $b_{55}$ 
of label $m+2$ 
are oriented from $w_5$ to $w_7$
(see Fig.~\ref{fig19}(a)),
both of $e_7'$ and $e_7''$ are oriented outward at $w_7$.
Moreover by Assumption~\ref{AssumeTerminal},
neither $e_7'$ nor $e_7''$ is a terminal edge.
Hence 
we have $w(\Gamma_{m+3}\cap{\rm Int}F_2)\ge1$
by IO-Calculation with respect to $\Gamma_{m+3}$ in $F_2$.
This contradicts (7).
Thus $e_7\not\subset F_2$.
Hence Claim~3 holds.
{\hfill {$\square$}\vspace{1.5em}}

Now we are ready to prove 
Proposition~\ref{PropOvalGammaM+1}.
There are two cases:
(i) $F_2\ni w_4$,
(ii) $F_2\not\ni w_4$.

{\bf Case (i).}
First we shall show that
for each $i=1,2$,
the intersection $e_i\cap \partial F_2$ is one point
by C-moves in a neighborhood of $F_2$
keeping $\Gamma_{m+1}\cup\Gamma_{m+2}\cup\Gamma_{m+3}$ 
fixed (see Fig.~\ref{fig11}(b)).

Since $w_4\in{\rm Int}F_2$ by the condition of Case (i),
we have $a_{44}\cup b_{44}\subset F_2$.
Let $E$ be the 2-angled disk of $\Gamma_{m+1}$ in $F_2$
with $\partial E=a_{44}\cup b_{44}$.
By (8),
we have 
\begin{enumerate}
\item[(9)] $w(\Gamma\cap({\rm Int}F_2 -E))=0$.
\end{enumerate}

Since for each $i=1,2$
the edge $e_i$ connects the white vertex $w_4$ in 
${\rm Int}F_2$ and the white vertex $w_1$ 
in the outside of $F_2$,
the intersection $e_i\cap F_2$
consists of an arc with $w_4$ and proper arcs of $F_2$.
Moreover $\partial E\subset \Gamma_{m+1}$ and (6) 
imply
\begin{enumerate}
\item[(10)] $e_i\cap F_2=e_i\cap Cl(F_2-E)$,
\item[(11)] $e_i\cap F_2$ consists of 
an arc with $w_4$ and $(F_2,b_{55})$-arcs of label $m$.
\end{enumerate}

By (7) and Claim~3,
there exists a terminal edge $e'$ of label $m+3$ at $w_7$.
Let $e''$ be the terminal edge of label $m+1$ at $w_4$,
and $N'$ a regular neighborhood of 
$E\cup e'\cup e''$ in $F_2$.
Then $\widetilde{F_2}=Cl(F_2-N')$ is a disk.
Let $\widetilde{b_{55}}=b_{55}\cap \widetilde{F_2}$.
By (9),
we can apply New Disk Lemma (Lemma~\ref{NewDiskLemma})
for the disk $\widetilde{F_2}$
so that the chart $\Gamma$ is 
$(\widetilde{F_2},\widetilde{b_{55}})$-arc free. 
Thus
there is no $(\widetilde{F_2},\widetilde{b_{55}})$-arc
of label $m$,
i.e. there is no $(F_2,b_{55})$-arcs of label $m$.
Thus by (10) and (11) for each $i=1,2$,
the intersection $e_i\cap F_2$ is an arc with $w_4$.
Therefore
the chart $\Gamma$ contains the pseudo chart as shown in 
Fig.~\ref{fig11}(b).
However by Lemma~\ref{NoMinimalOvalGammaM+1},
the chart $\Gamma$ is not minimal.
This is a contradiction.
Hence Case (i) does not occur.

{\bf Case (ii).}
Similarly we can show that 
$e_i\cap \partial F_2=\emptyset$ for each $i=1,2$
by C-moves in a neighborhood of $F_2$
keeping $\Gamma_{m+1}\cup\Gamma_{m+2}\cup\Gamma_{m+3}$ 
fixed.
Thus 
the chart $\Gamma$ contains the pseudo chart as shown in 
Fig.~\ref{fig11}(a).
However by Lemma~\ref{NoMinimalOvalGammaM+1},
the chart $\Gamma$ is not minimal.
This is a contradiction.
Hence Case (ii) does not occur.

Therefore the graph $\Gamma_{m+1}$ does not contain
an oval.
By Lemma~\ref{GammaM+1ToGammaM+2},
the graph $\Gamma_{m+2}$ does not contain
an oval.
We complete the proof of 
Proposition~\ref{PropOvalGammaM+1}.
\end{Proof}

\begin{figure}
\centerline{\includegraphics{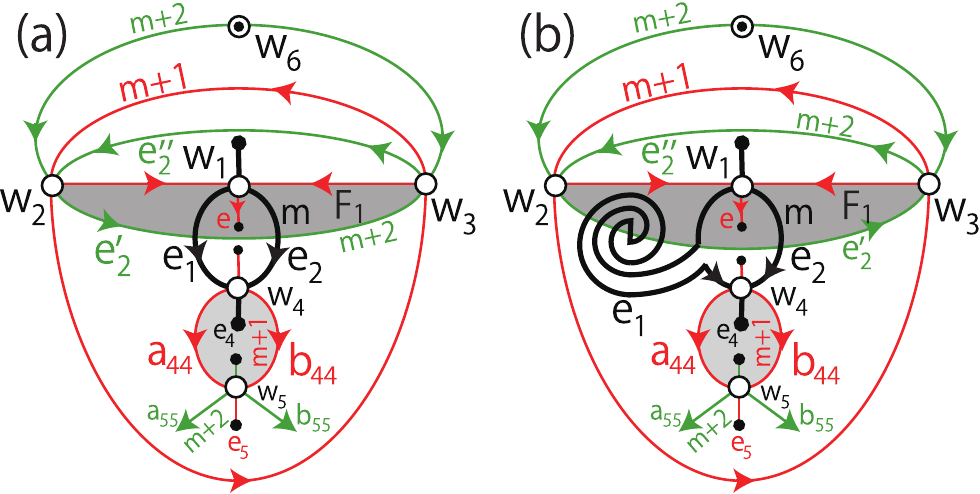}}
\caption{\LABEL{fig19}
The dark gray regions are the disk $F_1$.
The light gray regions are the disk $E$.
}
\end{figure}

By Lemma~\ref{OvalTypeGTypeH}
 and
Proposition~\ref{PropOvalGammaM+1},
we have the following corollary:

\begin{corollary}
\LABEL{TypeGTypeH}
If there exists a minimal chart $\Gamma$ of 
type $(m;2,3,2)$,
then each of $\Gamma_{m+1}$ and $\Gamma_{m+2}$
contains 
 one of two graphs as shown 
in  Fig.~\ref{fig05}.
\end{corollary}



\section{3-angled disks}
\LABEL{s:3angledDisk}

In this section,
we investigate a minimal chart $\Gamma$ 
of type $(m;2,3,2)$
such that $\Gamma_{m+1}$ contains 
the graph as shown in Fig.~\ref{fig05}(a).

Let $X$ be a set in a chart $\Gamma$. Let 
\begin{center}
$c(X)=${the number of crossings
on $X$}.
\end{center}

Let $D$ be a $k$-angled disk of $\Gamma_m$ for a minimal chart $\Gamma$.
The pair of integers 
$(w(\Gamma\cap{\rm Int}D),c(\Gamma\cap\partial D))$
is called the {\it local complexity 
with respect to $D$},
denoted by $\ell c(D;\Gamma)$.
Let
${\Bbb S}$ be the set of all minimal charts each of which can be moved from $\Gamma$ by C-moves in a regular neighborhood of $D$ keeping $\partial D$ fixed.
The chart $\Gamma$ is said to be 
{\it locally minimal
with respect to $D$}
if its local complexity
with respect to $D$
is minimal
among the charts in ${\Bbb S}$ with respect to 
the lexicographic order.

Let $\Gamma$ be a chart,
and $D$ a $k$-angled disk of $\Gamma_m$.
If any feeler of $D$ of label $m$ is a terminal edge,
then $D$ is called a {\it special} $k$-angled disk.

\begin{lemma}
{\rm (\cite[Theorem 1.2]{ChartAppIII})}
\LABEL{Theorem3AngledDisk}
Let $\Gamma$ be a minimal chart.
Let $D$ be a special $3$-angled disk of $\Gamma_m$
such that $\Gamma$ is locally minimal
with respect to $D$.
If $w(\Gamma\cap {\rm Int}D)\le1$,
then a regular neighborhood of $D$ contains one of the RO-families of the eight pseudo charts as shown in Fig.~\ref{fig20}.
\end{lemma}

\begin{figure}[htb]
\centerline{\includegraphics{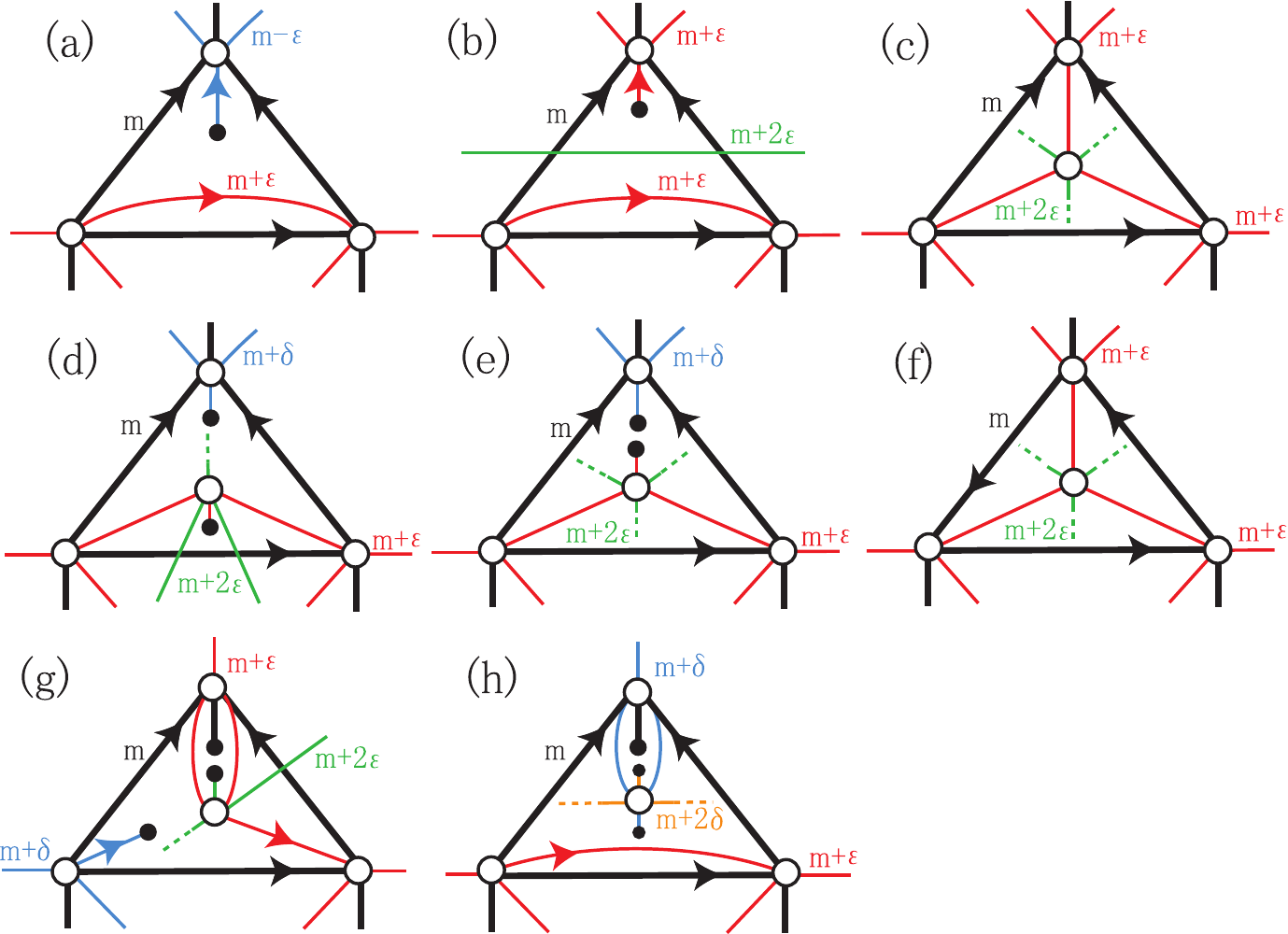}}
\caption{\LABEL{fig20} 
The 3-angled disks (g) and (h) have one feeler,
 the others do not have any feelers.
$m$ is a label, $\varepsilon,\delta\in\{+1,-1\}$.}
\end{figure}

\begin{lemma}
\LABEL{LemmaTypeGDiskD1}
Let $\Gamma$ be a minimal chart of type $(m;2,3,2)$.
Suppose that $\Gamma_{m+1}$ contains the graph as shown
in Fig.~\ref{fig05}$($a$)$.
Let $D_1$ be a $3$-angled disk of $\Gamma_{m+1}$
with $w(\Gamma_{m+1}\cap{\rm Int}D_1)=0$.
If $w(\Gamma\cap{\rm Int}D_1)\le 1$
and if $\Gamma$ is locally minimal with respect to $D_1$,
then $\Gamma$ contains one of the RO-families of 
the three pseudo charts as shown in 
Fig.~\ref{fig21}.
\end{lemma}

\begin{figure}[htb]
\centerline{\includegraphics{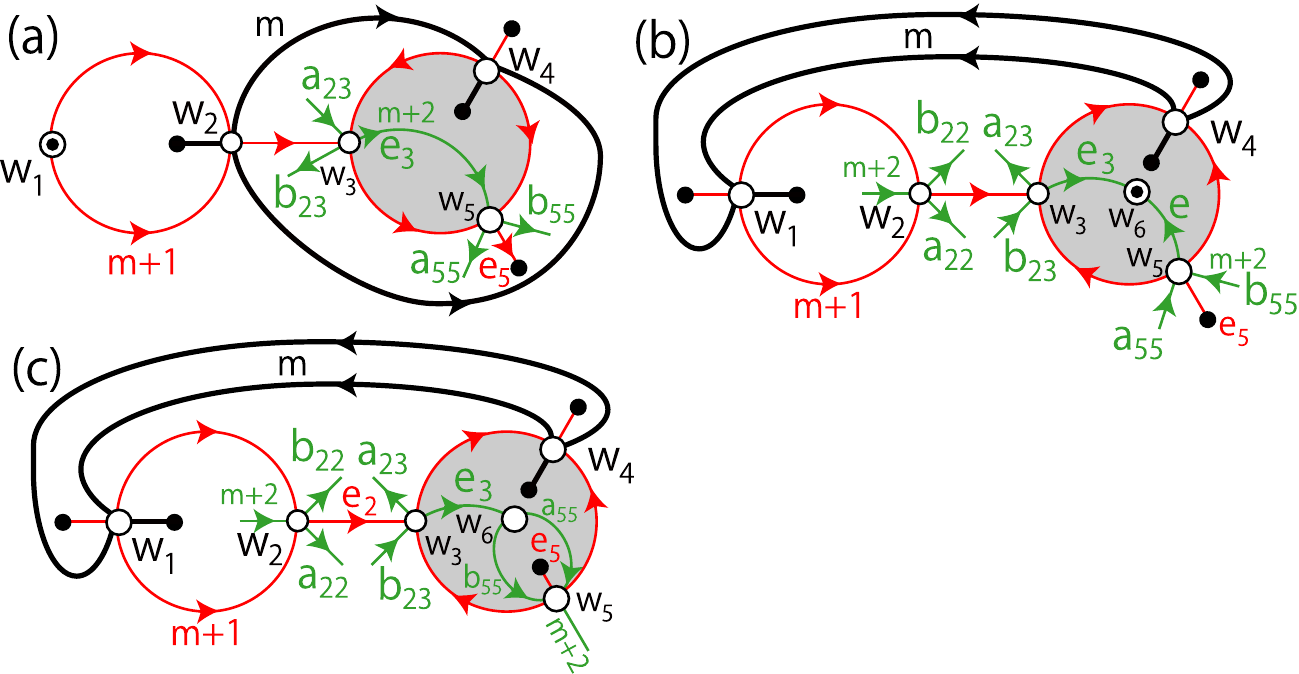}}
\caption{\LABEL{fig21} 
The gray regions are the disk $D_1$.
}
\end{figure}

\begin{Proof}
By Lemma~\ref{LemmaTypeG},
the chart $\Gamma$ contains one of the RO-families of 
the two pseudo charts as shown in 
Fig.~\ref{fig08}(a),(b).
We use the notations as shown in 
Fig.~\ref{fig08}(a),(b),
where
\begin{enumerate}
\item[(1)] $e_3$ is an internal edge 
(possibly terminal edge) of label $m+2$ 
oriented outward at $w_3$,
\item[(2)] $w_3\in\Gamma_{m+2}$,
$w_4\in\Gamma_m$.
\end{enumerate}

We shall show $w_5\in\Gamma_{m+2}$.
Since $w_5\in\Gamma_{m+1}$,
we have $w_5\in\Gamma_m$ or $w_5\in\Gamma_{m+2}$.
Suppose $w_5\in\Gamma_m$.
Since one of $w_1,w_2$ is in $\Gamma_m$, by (2)
the graph $\Gamma_m$ contains at least three
white vertices.
This contradicts the fact that $\Gamma$ 
is of type $(m;2,3,2)$.
Hence we have
$w_5\in\Gamma_{m+2}$.

Since $w_4\in\Gamma_m$ and $w_3,w_5\in\Gamma_{m+2}$
by (2), 
a regular neighborhood $N(D_1)$ of $D_1$
contains one of the RO-families of 
the five pseudo charts as shown in 
Fig.~\ref{fig20}(a),(d),(e),(g),(h)
by Lemma~\ref{Theorem3AngledDisk}.

Looking at Fig.~\ref{fig08}(a),(b),
we see that
the terminal edge of label $m+1$ at $w_4$
is not contained in $D_1$.
Thus $N(D_1)$ does not contain one of the RO-family of 
the pseudo chart as shown in 
Fig.~\ref{fig20}(h).

Let $e_5$ be the terminal edge of label $m+1$ at $w_5$.

If $N(D_1)$ contains one of the RO-family of 
the pseudo chart as shown in 
Fig.~\ref{fig20}(a),
then $e_5\not\subset D_1$ and $e_3\ni w_5$.
Thus by (1),
the edge $e_3$ is oriented inward at $w_5$.
Hence $\Gamma$ contains the pseudo chart as shown
in Fig.~\ref{fig08}(b).
Thus $\Gamma$ contains the pseudo chart as shown
in Fig.~\ref{fig21}(a).

If $N(D_1)$ contains one of the RO-families of 
the two pseudo charts as shown in 
Fig.~\ref{fig20}(d),(e),
then $e_5\not\subset D_1$. 
Moreover 
there exists only one internal edge $e$ of label $m+2$
at $w_5$ in $D_1$ such that $e\cap e_3$
is a BW-vertex with respect to $\Gamma_{m+2}$.
Thus by (1) and Lemma~\ref{OriBWvertex},
the edge $e$ is oriented outward at $w_5$.
Hence $\Gamma$ contains the pseudo chart as shown
in Fig.~\ref{fig08}(a).
Thus $\Gamma$ contains the pseudo chart as shown
in Fig.~\ref{fig21}(b).

If $N(D_1)$ contains one of the RO-family of 
the pseudo chart as shown in 
Fig.~\ref{fig20}(g),
then $e_5\subset D_1$ and 
there exist two internal edges $a_{55},b_{55}$ of 
label $m+2$ at $w_5$ in $D_1$
such that
$e_3\cap a_{55}\cap b_{55}$ is a white vertex.
Thus by (1) and Assumption~\ref{AssumeTerminal},
the edges $a_{55},b_{55}$ are oriented inward at $w_5$.
Hence $\Gamma$ contains the pseudo chart as shown
in Fig.~\ref{fig08}(a).
Thus $\Gamma$ contains the pseudo chart as shown
in Fig.~\ref{fig21}(c).
\end{Proof}



\section{Case of the graph as shown in Fig.~\ref{fig05}(a)}
\label{s:CaseG}

In this section
we shall show that
neither $\Gamma_{m+1}$ nor $\Gamma_{m+2}$ contains
the graph as shown in Fig.~\ref{fig05}(a)
for any minimal chart $\Gamma$ of type $(m;2,3,2)$.

Let $\Gamma$ be a chart,
and $m$ a label of $\Gamma$. 
A {\it loop} is a simple closed curve in $\Gamma_m$ with exactly one white vertex
(possibly with crossings).

\begin{lemma}$(${\rm \cite[Theorem 1.1]{ChartAppIV}}$)$
\LABEL{LemmaNoLoop}
There is no loop in any minimal chart with exactly seven white vertices.
\end{lemma}

\begin{lemma}
\LABEL{LemmaTypeGTypeA}
Let $\Gamma$ be a chart of type $(m;2,3,2)$.
If $\Gamma$ contains the pseudo chart as shown in 
Fig.~\ref{fig08}$($a$)$,
then $\Gamma$ is not minimal.
\end{lemma}

\begin{Proof}
Suppose that $\Gamma$ is minimal.
We use the notations as shown in 
Fig.~\ref{fig08}(a),
where
\begin{enumerate}
\item[(1)] $w_2\in\Gamma_{m+1}\cap\Gamma_{m+2}$,
\item[(2)] $a_{22},b_{22},a_{23}$ are internal edges
(possibly terminal edges) of label $m+2$ 
oriented outward at $w_2,w_2,w_3$
but not middle at $w_2,w_2,w_3$, respectively,
\item[(3)] $b_{23}$ is an internal edge
(possibly terminal edge) of label $m+2$ 
middle at $w_3$.
\end{enumerate}
By (2) and Assumption~\ref{AssumeTerminal},
\begin{enumerate}
\item[(4)] none of $a_{22},b_{22},a_{23}$ 
are terminal edges.
\end{enumerate}

Let $D_1$ be the 3-angled disk of $\Gamma_{m+1}$ with
$\partial D_1\ni w_3,w_4,w_5$ and
$D_1\supset e_3$ (see Fig.~\ref{fig08}(a)).
We can assume that
\begin{enumerate}
\item[(5)] $\Gamma$ is locally minimal 
with respect to $D_1$.
\end{enumerate}

{\bf Claim 1.} $w(\Gamma\cap{\rm Int}D_1)\ge1$.

{\it Proof of Claim~$1$.}
Suppose that $w(\Gamma\cap{\rm Int}D_1)=0$.
Since $\Gamma_{m+1}$ contains the graph as shown
in Fig.~\ref{fig05}(a),
by (5) and Lemma~\ref{LemmaTypeGDiskD1}
the chart $\Gamma$ contains one of the RO-family
of the pseudo chart as shown in 
Fig.~\ref{fig21}(a).
Thus $w_2\in\Gamma_m$.
This contradicts (1). 
Hence we have $w(\Gamma\cap{\rm Int}D_1)\ge1$.
Thus Claim~1 holds. {\hfill {$\square$}\vspace{1.5em}}

Let $D_2$ be the 2-angled disk of $\Gamma_{m+1}$
with $\partial D_2\ni w_1,w_2$ and 
$D_1\cap D_2=\emptyset$ 
(see Fig.~\ref{fig08}(a)).

{\bf Claim~2.}
$w(\Gamma\cap(S^2-(D_1\cup D_2)))\ge1$.

{\it Proof of Claim~$2$.}
Suppose that $w(\Gamma\cap(S^2-(D_1\cup D_2)))=0$.
Since $a_{22},b_{22},a_{23}$ are oriented outward 
at $w_2,w_2,w_3$ respectively by (2),
and since none of $a_{22},b_{22},a_{23}$ are
terminal edges by (4),
we have $a_{22}=b_{23}$ and
$b_{22}\cap a_{23}\ni w_5$
(see Fig.~\ref{fig22}(a)).
Thus the three white vertices $w_2,w_3,w_5$
are contained in the simple closed curve 
$C=a_{22}\cup b_{22}\cup a_{23}$ 
of $\Gamma_{m+2}$ with $w(C)=3$.
Hence by Corollary~\ref{TypeGTypeH},
the graph $\Gamma_{m+2}$ contains 
the graph as shown in Fig.~\ref{fig05}(a).
Thus by Lemma~\ref{LemmaTypeGTypeHGammaM+2},
the simple closed curve $C$ contains one white vertex
in $\Gamma_{m+3}$
(see Fig.~\ref{fig09}(a),(b)).
However all of the white vertices $w_2,w_3,w_5$ in $C$
are contained in $\Gamma_{m+1}\cap\Gamma_{m+2}$.
This is a contradiction.
Hence $w(\Gamma\cap(S^2-(D_1\cup D_2)))\ge1$.
Thus Claim~2 holds. {\hfill {$\square$}\vspace{1.5em}}

{\bf Claim~3.}
$w(\Gamma\cap {\rm Int}D_1)=1$ and 
$w(\Gamma\cap(S^2-(D_1\cup D_2)))=1$.

{\it Proof of Claim~$3$.}
Since $D_2$ is a 2-angled disk,
we have $w(\Gamma\cap D_2)\ge2$.
Since $D_1$ is a 3-angled disk,
we have $w(\Gamma\cap\partial D_1)=3$.
Thus 
$$w(\Gamma\cap D_1)=w(\Gamma\cap\partial D_1)+w(\Gamma\cap{\rm Int} D_1)=3+w(\Gamma\cap{\rm Int} D_1).$$
Hence
$$\begin{array}{rl}
7=w(\Gamma) & =w(\Gamma\cap D_1)+w(\Gamma\cap D_2)
+w(\Gamma\cap(S^2-(D_1\cup D_2)))\\
& \ge 3+w(\Gamma\cap {\rm Int}D_1)+2
+w(\Gamma\cap(S^2-(D_1\cup D_2))).\\
\end{array}$$
Thus
$$2\ge w(\Gamma\cap {\rm Int}D_1)
+w(\Gamma\cap(S^2-(D_1\cup D_2))).$$
Hence Claim~1 and Claim~2 imply Claim~3. {\hfill {$\square$}\vspace{1.5em}}

By Claim~3,
there exists a white vertex in ${\rm Int}D_1$, 
say $w_6$, and
there exists a white vertex in
$S^2-(D_1\cup D_2)$, say $w_7$.
Since $\Gamma$ is of type $(m;2,3,2)$
and since $w_1,w_2,w_3,w_4,w_5\in \Gamma_{m+1}$,
we have  
\begin{enumerate}
\item[(6)] the graph $\Gamma_{m+3}$ contains 
only two white vertices $w_6,w_7$.
\end{enumerate}

Since $w(\Gamma\cap {\rm Int}D_1)=1$ by Claim~3,
the chart $\Gamma$ contains one of the RO-families
of the two pseudo charts as shown in 
Fig.~\ref{fig21}(b),(c)
by (5) and Lemma~\ref{LemmaTypeGDiskD1}.
Hence there are two cases.

{\bf Case (i).}
Suppose that 
the chart $\Gamma$ contains one of the RO-family
of the pseudo chart as shown in 
Fig.~\ref{fig21}(b).
We use the notations as shown in 
Fig.~\ref{fig21}(b),
where
\begin{enumerate}
\item[(7)] both internal edges $e_3,e$ of label $m+2$
are oriented inward at $w_6$, and
\end{enumerate}
$a_{55},b_{55}$ are internal edges
(possibly terminal edges) of label $m+2$
 at $w_5$ not middle at $w_5$.
Thus by Assumption~\ref{AssumeTerminal},
\begin{enumerate}
\item[(8)] neither $a_{55}$ nor $b_{55}$ is 
a terminal edge.
\end{enumerate}

By Corollary~\ref{TypeGTypeH},
the graph $\Gamma_{m+2}$ contains
one of the two graphs as shown in
Fig.~\ref{fig05}.
Hence there are two cases.

If $\Gamma_{m+2}$ contains the graph as shown in 
Fig.~\ref{fig05}(a),
then the white vertex $w_6$ must contain
the simple closed curve $C$
of label $m+2$ with $w(C)=3$.
Thus by (8), we have
$b_{55}=a_{23}$ and $C=b_{55}\cup e_3\cup e$.
Hence the edge $b_{23}$ is a terminal edge,
 $a_{22}\cap b_{22}\cap a_{55}\ni w_7$ 
(see Fig.~\ref{fig22}(b)).

Since both of $a_{22}$ and $b_{22}$ of label $m+2$
are oriented from $w_2$ to $w_7$ by (2),
there exist two internal edges $e_7',e_7''$
(possibly terminal edges) of label $m+3$
oriented outward at $w_7$
but not middle at $w_7$.
Thus by Assumption~\ref{AssumeTerminal},
neither $e_7'$ nor $e_7''$ is a terminal edge.
Hence by Lemma~\ref{LemmaNoLoop} and (6),
we have $e_7'\cap e_7''\ni w_6$.
Thus the two edges $e_7',e_7''$ of label $m+3$
are oriented inward at $w_6$.
However by (7)
there exist four edges oriented inward at $w_6$.
This contradicts Condition (iii) of
the definition of a chart.
Hence $\Gamma_{m+2}$ does not contain 
the graph as shown in Fig.~\ref{fig05}(a).

If $\Gamma_{m+2}$ contains 
the graph as shown in Fig.~\ref{fig05}(b),
then the BW-vertex $w_6$ with respect to $\Gamma_{m+2}$
is not contained 
in any simple closed curve of label $m+2$
(Fig.~\ref{fig21}(b)).
Hence both of $a_{23}\cup b_{23}$ and $a_{55}\cup b_{55}$
are simple closed curves.
Thus by Lemma~\ref{LemmaTypeGTypeHGammaM+2},
any simple closed curve of label $m+2$
does not contain middle arcs at any white vertex
(see Fig.~\ref{fig09}(c),(d)).
Hence the edge $b_{23}$ is not middle at $w_3$.
This contradicts (3).
Thus $\Gamma_{m+2}$ does not contain 
the graph as shown in Fig.~\ref{fig05}(b).

Hence neither two cases occur.
Thus Case (i) does not occur.

{\bf Case (ii).}
Suppose that 
the chart $\Gamma$ contains one of the RO-family
of the pseudo chart as shown in 
Fig.~\ref{fig21}(c).
We use the notations as shown in 
Fig.~\ref{fig21}(c),
where
\begin{enumerate}
\item[(9)] $a_{55}\cup b_{55}$ is a simple closed curve.
\end{enumerate}

By Corollary~\ref{TypeGTypeH},
the graph $\Gamma_{m+2}$ contains one of 
the two graphs as shown in Fig.~\ref{fig05}.
Hence there are two cases.

If $\Gamma_{m+2}$ contains the graph as shown 
in Fig.~\ref{fig05}(a),
then there exists a simple closed curve $C$ 
of label $m+2$ with $w(C)=3$.
Thus by (9),
$w_5\not\in C$ and $w_6\not\in C$.
Hence $C\ni w_2,w_3,w_7$ and
$C=a_{22}\cup b_{22}\cup a_{23}$.
Moreover there exist terminal edges 
of label $m+2$ at $w_2,w_7$ respectively.

Let $D_3$ be the 3-angled disk of $\Gamma_{m+2}$
with $\partial D_3=C$ 
which contains 
the internal edge $e_2$ of label $m+1$
connecting $w_2,w_3$.
By (6) and Corollary of Triangle Lemma
(Corollary~\ref{CorollaryTriangleLemma}),
we have $w(\Gamma\cap {\rm Int}D_3)\ge1$.
Thus there exists a white vertex
in ${\rm Int}D_3$ different from $w_7$.
Thus the condition
${\rm Int}D_3\subset S^2-(D_1\cup D_2)$
implies $w(\Gamma\cap(S^2-(D_1\cup D_2)))\ge2$.
This contradicts Claim~3.
Hence $\Gamma_{m+2}$ does not contain 
the graph as shown in Fig.~\ref{fig05}(a).

If $\Gamma_{m+2}$ contains the graph as shown 
in Fig.~\ref{fig05}(b),
then by (9)
one of the edges $a_{23},b_{23}$ is a terminal edge.
Since $a_{23}$ is not a terminal edge by (4),
the edge $b_{23}$ is a terminal edge.
Hence $a_{22}\cap b_{22}\cap a_{23}\ni w_7$.
However by (2),
all of the three edges $a_{22}, b_{22}, a_{23}$
are edges of label $m+2$ oriented inward at $w_7$.
This contradicts Condition (iii)
of the definition of a chart.
Hence $\Gamma_{m+2}$ does not contain 
the graph as shown in Fig.~\ref{fig05}(b).

Thus neither two cases occur.
Hence Case (ii) does not occur.

Therefore neither
Case (i) nor Case (ii) occurs.
Hence $\Gamma$ is not minimal.
We complete the proof of Lemma~\ref{LemmaTypeGTypeA}.
\end{Proof}

\begin{figure}[htb]
\centerline{\includegraphics{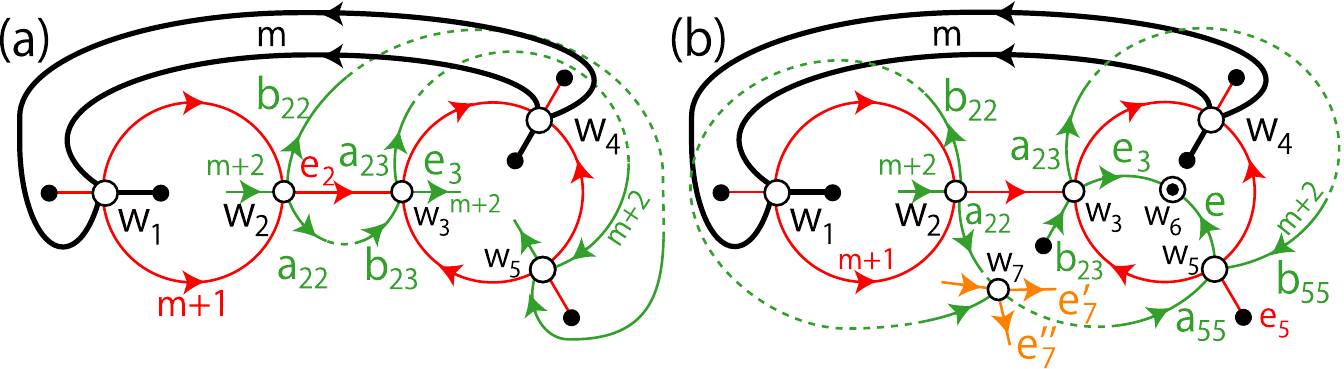}}
\caption{\LABEL{fig22}
(a) There is a simple closed curve $a_{22}\cup b_{22}\cup a_{23}$ in $\Gamma_{m+2}$.
(b) The graph $\Gamma_{m+2}$ contains the graph as shown in
Fig.~\ref{fig05}(a).
}
\end{figure}

\begin{proposition}
\LABEL{LemmaTypeGTypeB}
Let $\Gamma$ be a minimal chart of type $(m;2,3,2)$.
Then neither of $\Gamma_{m+1}$ nor $\Gamma_{m+2}$
contains the graph as shown in 
Fig.~\ref{fig05}$($a$)$.
\end{proposition}

\begin{Proof}
Suppose that $\Gamma_{m+1}$ contains
the graph as shown in Fig.~\ref{fig05}$($a$)$.
By Lemma~\ref{LemmaTypeG} and 
Lemma~\ref{LemmaTypeGTypeA},
the chart $\Gamma$ contains one of the RO-family of
the pseudo chart 
as shown in Fig.~\ref{fig08}(b).
We use the notations as shown in 
Fig.~\ref{fig08}(b), where
\begin{enumerate}
\item[(1)] $w_2\in\Gamma_m$.
\end{enumerate}
Moreover,
$b_{23}$ is an internal edge 
(possibly terminal edge) of label $m+2$ at $w_3$
but not middle at $w_3$.
By Assumption~\ref{AssumeTerminal},
\begin{enumerate}
\item[(2)] $b_{23}$ is not a terminal edge. 
\end{enumerate}

Let $D_1$ be the 3-angled disk of $\Gamma_{m+1}$
with $\partial D_1\ni w_3,w_4,w_5$ and
$D_1\supset e_3$ 
(see Fig.~\ref{fig08}(b)).

{\bf Claim~1.}
$w(\Gamma\cap{\rm Int}D_1)\ge2$.

{\it Proof of Claim~$1$.}
Suppose that $w(\Gamma\cap{\rm Int}D_1)\le1$.
We can assume that
$\Gamma$ is locally minimal 
with respect to $D_1$. 
Thus (1) 
and Lemma~\ref{LemmaTypeGDiskD1} imply that
the chart $\Gamma$ contains
one of the RO-family of the pseudo chart
as shown in Fig.~\ref{fig21}(a).
Since a neighborhood of $D_1$ 
contains the pseudo chart as shown in
Fig.~\ref{fig13}(a),
we can apply Triangle Lemma (Lemma~\ref{LemmaTriangle})
so that $w(\Gamma\cap{\rm Int}D_1)\ge1$.
Hence $w(\Gamma\cap{\rm Int}D_1)=1$.

Let $w_6$ be the white vertex in ${\rm Int}D_1$.
Since $w_1,w_2,w_3,w_4,w_5\in\Gamma_{m+1}$ and
since $\Gamma$ is of type $(m;2,3,2)$,
we have $w_6\in\Gamma_{m+2}\cap\Gamma_{m+3}$.
Hence there exists a connected component of $\Gamma_{m+2}$
with only one white vertex $w_6$.
This contradicts Lemma~\ref{LemmaWithTerminal3}.
Thus $w(\Gamma\cap{\rm Int}D_1)\ge2$.
Hence Claim~1 holds. {\hfill {$\square$}\vspace{1.5em}}

Let $D_2$ be the 2-angled disk of $\Gamma_{m+1}$
with $\partial D_2\ni w_1,w_2$
and $D_1\cap D_2=\emptyset$ 
(see Fig.~\ref{fig08}(b)).

{\bf Claim~2.}
$w(\Gamma\cap {\rm Int} D_2)=0$
and
$w(\Gamma\cap (S^2-(D_1\cup D_2)) )=0$.

{\it Proof of Claim~$2$.}
Since $D_1$ is a 3-angled disk, by Claim~1
we have 
$$w(\Gamma\cap D_1)=w(\Gamma\cap \partial D_1)+w(\Gamma\cap {\rm Int} D_1)\ge 3+2=5.$$
Since $D_2$ is a 2-angled disk,
we have 
$$w(\Gamma\cap D_2)=w(\Gamma\cap \partial D_2)+w(\Gamma\cap {\rm Int} D_2)=2+w(\Gamma\cap {\rm Int} D_2).$$
Hence 
$$\begin{array}{rl}
7=w(\Gamma) & =w(\Gamma\cap D_1)+w(\Gamma\cap D_2)
+w(\Gamma\cap(S^2-(D_1\cup D_2)))\\
& \ge 5+2+w(\Gamma\cap {\rm Int} D_2)+w(\Gamma\cap(S^2-(D_1\cup D_2))).\\
\end{array}$$
Thus $0\ge w(\Gamma\cap {\rm Int} D_2)+w(\Gamma\cap(S^2-(D_1\cup D_2))).$
Therefore Claim~2 holds. {\hfill {$\square$}\vspace{1.5em}}

Since $w(\Gamma\cap {\rm Int} D_2)=0$ by Claim~2,
a neighborhood of $D_2$ contains the pseudo chart
as shown in Fig.~\ref{fig10}(a)
by Lemma~\ref{Theorem2AngledDisk}.
Let $e_1$ be the terminal edge of label $m+1$ at $w_1$.
Moreover,
we have the following claim:

{\bf Claim~3.} $e_1\not\subset D_2$
(see Fig.~\ref{fig23}).

Let $a_{11},b_{11}$ be internal edges 
(possibly terminal edges) of label $m+2$
such that $a_{11},e_1,b_{11}$ are oriented inward at $w_1$
and lie anticlockwise around $w_1$.
Then by Assumption~\ref{AssumeTerminal}
we have
\begin{enumerate}
\item[(3)] neither $a_{11}$ nor $b_{11}$ 
is a terminal edge. 
\end{enumerate}

Let $e_5$ be the terminal edge of label $m+1$ at $w_5$.

{\bf Claim~4.} $e_5\not\subset D_1$.

{\it Proof of Claim~$4$.}
Suppose $e_5\subset D_1$.
Then there exists only one edge $e$
of label $m+2$ at $w_5$ in $Cl(S^2-(D_1\cup D_2))$
(see Fig.~\ref{fig23}(a)).
Moreover
\begin{enumerate}
\item[(4)] the edge $e$ is oriented inward at $w_5$. 
\end{enumerate}
Since $a_{11},b_{11},a_{23}$ are oriented inward at
$w_1,w_1,w_3$ respectively 
(see Fig.~\ref{fig23}(a)),
and since neither $a_{11}$ nor $b_{11}$
is a terminal edge by (3),
we have $w(\Gamma\cap(S^2-(D_1\cup D_2)))\ge1$
by IO-Calculation with respect to $\Gamma_{m+2}$
in $Cl(S^2-(D_1\cup D_2))$.
This contradicts Claim~2.
Hence $e_5\not\subset D_1$.
Thus Claim~4 holds. {\hfill {$\square$}\vspace{1.5em}}

Let $a_{55},b_{55}$ be internal edges
(possibly terminal edges) of label $m+2$ at $w_5$
such that $a_{55},e_{5},b_{55}$ are
consecutive three edges lying anticlockwise 
around $w_5$.
Then by Assumption~\ref{AssumeTerminal},
\begin{enumerate}
\item[(5)] neither $a_{55}$ nor $b_{55}$ 
is a terminal edge. 
\end{enumerate}

Now we are ready to prove 
Proposition~\ref{LemmaTypeGTypeB}.
By Claim~3 and Claim~4,
the chart $\Gamma$ contains the pseudo chart
as shown in Fig.~\ref{fig23}(b).
Since $b_{23},a_{55},b_{55}$ are oriented outward at 
$w_3,w_5,w_5$ respectively
(see Fig.~\ref{fig23}(b)),
and since none of $b_{23},a_{55},b_{55}$
are terminal edges by (2) and (5),
the condition $w(\Gamma\cap(S^2-(D_1\cup D_2)))=0$
implies 
$b_{23}=b_{11}$, $a_{55}=a_{11}$
and $b_{55}=a_{23}$.
Thus 
the simple closed curve 
$C=b_{23}\cup a_{55}\cup b_{55}$
of label ${m+2}$ 
contains only three white vertices
$w_1,w_3,w_5$.
Hence by Corollary~\ref{TypeGTypeH},
the graph $\Gamma_{m+2}$ contains
the graph as shown in Fig.~\ref{fig05}(a).
Thus by Lemma~\ref{LemmaTypeGTypeHGammaM+2},
the simple closed curve $C$ contains one white vertex
of $\Gamma_{m+3}$
(see Fig.~\ref{fig09}(a),(b)).
However all of the white vertices $w_1,w_3,w_5$
in $C$ are contained in $\Gamma_{m+1}\cap\Gamma_{m+2}$.
This is a contradiction.
Therefore $\Gamma_{m+1}$ does not contain
the graph as shown in Fig.~\ref{fig05}(a).

By Lemma~\ref{GammaM+1ToGammaM+2},
we can show that $\Gamma_{m+2}$ does not contain
the graph as shown in Fig.~\ref{fig05}(a).
Therefore we complete the proof of
Proposition~\ref{LemmaTypeGTypeB}.
\end{Proof}

\begin{figure}[htb]
\centerline{\includegraphics{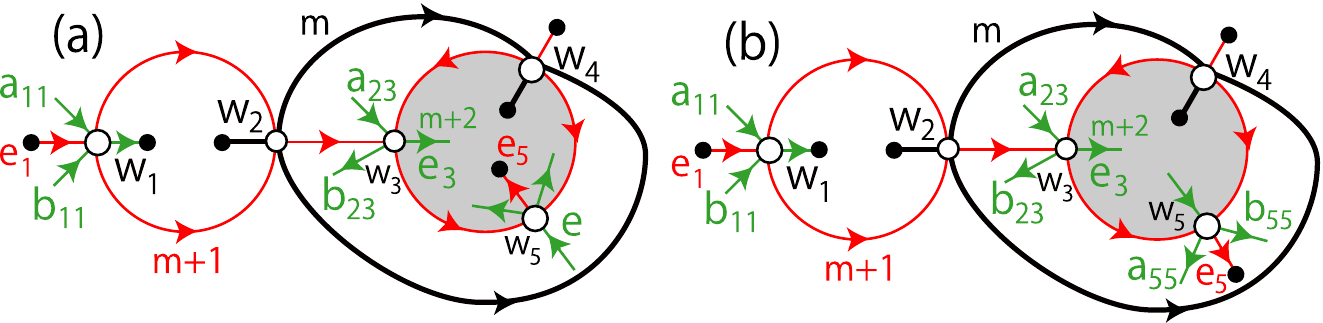}}
\caption{\LABEL{fig23}
The gray regions are the disk $D_1$.
(a) $e_5\subset D_1$.
(b) $e_5\not\subset D_1$.
}
\end{figure}


\section{Shifting Lemma}
\LABEL{s:ShiftingLemma}

In this section
we give an example of a non minimal chart $\Gamma$
of type $(m;2,3,2)$
such that $\Gamma_{m+1}$ contains 
the graph as shown in Fig.~\ref{fig05}(b).

Let $\Gamma$ be a chart.
Let $\alpha$ be an arc in an edge of $\Gamma_m$, 
and 
$w$ a white vertex with $w\not\in\alpha$. 
Suppose that there exists an arc $\beta$ in $\Gamma$ 
such that
its end points are the white vertex $w$ 
and an interior point $p$ of the arc $\alpha$.
Then we say that 
{\it the white vertex $w$ connects with the point $p$ of $\alpha$ by the arc $\beta$}.

Let $\alpha$ be a simple arc,  
and $p,q$ points in $\alpha$. 
We denote by $\alpha[p,q]$ the subarc of $\alpha$ whose endpoints are $p$ and $q$.

\begin{lemma}
{\rm (\cite[Lemma 4.2]{ChartApp1})}
$($Shifting Lemma$)$
\LABEL{Shift} 
Let $\Gamma$ be a chart and
$\alpha$ an arc in an edge of $\Gamma_m$.
Let
$w$ be a white vertex of $\Gamma_k \cap\Gamma_{h}$ where $h=k+\varepsilon ,\varepsilon\in\{+1,-1\}$.
Suppose that 
the white vertex $w$ connects with a point $r$ 
of the arc $\alpha$ by an arc in an edge $e$ of $\Gamma_k$. 
Suppose that
one of the following two conditions is satisfied: 
\begin{enumerate}
\item[{\rm (1)}] $h>k>m$.
\item[{\rm (2)}] $h<k<m$.
\end{enumerate}
Then for any neighborhood $V$ of the arc $e[w,r]$
we can shift the white vertex $w$
to $e-e[w,r]$ 
along the edge $e$
by C-I-R2 moves,
C-I-R3 moves and C-I-R4 moves in $V$ 
keeping $\displaystyle 
\bigcup_{i<0}\Gamma_{k+i\varepsilon}
$ 
fixed $($see Fig.~\ref{fig24}$)$.
\end{lemma}

\begin{figure}[htb]
\centerline{\includegraphics{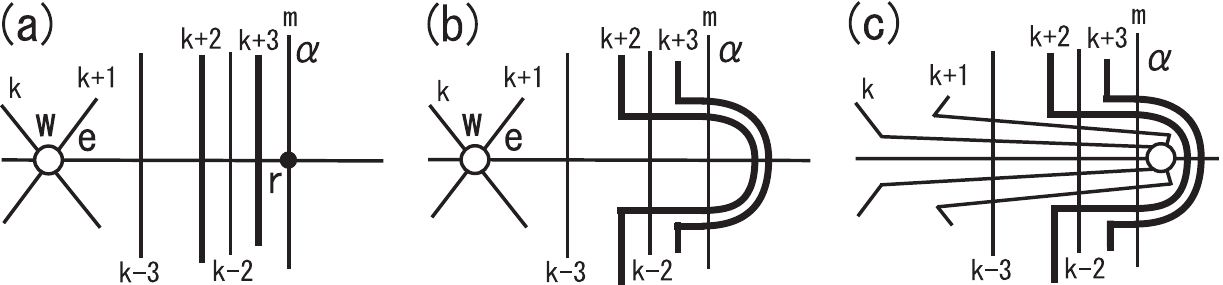}}
\caption{\LABEL{fig24}
$k>m$ and $\varepsilon=+1$.}
\end{figure}

\begin{lemma}
\LABEL{NoMinimalTypeH}
Let $\Gamma$ be a chart of type $(m;2,3,2)$
containing the pseudo chart as shown in 
Fig.~\ref{fig08}$($d$)$.
Let $e_2',e_2'',e_3',e_3'',e_5',e_5''$
be internal edges of label $m+2$
such that
$e_2',e_2''$ are oriented outward at the white vertex $w_2$,
$e_3',e_3''$ are oriented outward at the white vertex $w_3$, and
$e_5',e_5''$ are oriented inward at the white vertex $w_5$.
If $e_2'\ni w_3$,
$e_3'\cap e_3''$ contains a white vertex
different from $w_3$, and
$e_2''\cap e_5'\cap e_5''$ contains a white vertex
$($see Fig.~\ref{fig25}$($a$))$,
then $\Gamma$ is not a minimal chart.
\end{lemma}

\begin{figure}
\centerline{\includegraphics{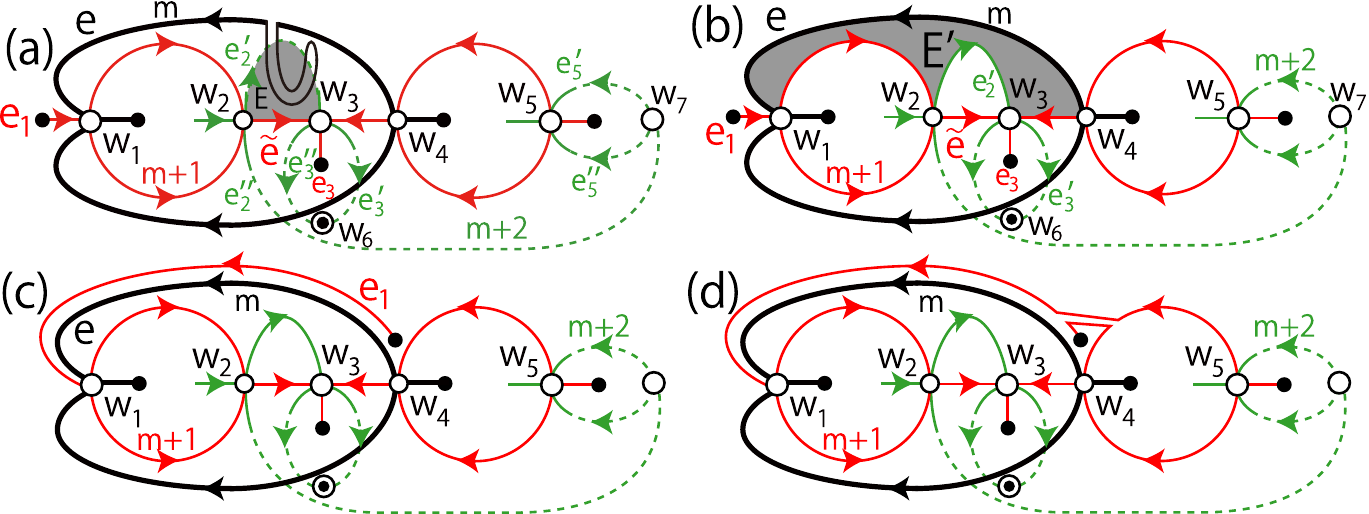}}
\caption{\LABEL{fig25}
(a) The gray region is the disk $E$.
(b) The gray region is the disk $E'$.
(c), (d) We apply C-II moves and a C-I-M2 move for the terminal edge $e_1$.}
\end{figure}

\begin{Proof}
Let $w_6,w_7$ be the white vertices different from
$w_2,w_3,w_5$ with
$w_6\in e_3'\cap e_3''$ and
$w_7\in e_2''\cap e_5'\cap e_5''$.
Let $\widetilde e$ be the internal edge of 
label $m+1$ with $\widetilde e \ni w_2,w_3$.
Let $e_3$ be the terminal edge of label $m+1$ at $w_3$.
Let $E$ be the disk with
$\partial E=e_2'\cup \widetilde e$
and $e_3\cap {\rm Int}E=\emptyset$.
Then
\begin{enumerate}
\item[(1)] none of the five white vertices $w_1,w_2,\cdots, w_5$ are contained in ${\rm Int}E$, 
\item[(2)] $(e_2''\cup e_3')\cap {\rm Int}E=\emptyset$.
\end{enumerate}

{\bf Claim 1.}
$w(\Gamma\cap{\rm Int}E)=0$.

{\it Proof of Claim $1$.}
By (1),
it suffices to prove that
neither $w_6$ nor $w_7$ is contained in ${\rm Int}E$.
Since the white vertices $w_6$ and $w_7$
are contained in the edges $e_3'$ and $e_2''$
respectively,
we have $w_6\not\in {\rm Int}E$
and $w_7\not\in {\rm Int}E$
by (2).
Thus Claim~1 holds.
{\hfill {$\square$}\vspace{1.5em}}

By Claim~1,
we can apply New Disk Lemma 
(Lemma~\ref{NewDiskLemma}) 
for the disk $E$ so that
we have the following claim:

{\bf Claim 2.}
We can assume that $\Gamma$ is $(E,e_2')$-arc free
by C-moves in a neighborhood of $e_2'$
keeping $\Gamma_{m+1}\cup\Gamma_{m+2}$ fixed.

Let $G,G_{m+1}$ be the connected components
of $\Gamma_m,\Gamma_{m+1}$ respectively
such that $w(G)=2$ and $w(G_{m+1})=5$.
Let $F$ be the closure of the connected component of
$S^2-(G\cup G_{m+1})$
with $\partial F\ni w_1,w_2,w_3,w_4$
and $e_3\cap{\rm Int}F=\emptyset$.
Then 
\begin{enumerate}
\item[(3)] none of the five white vertices
$w_1,w_2,\cdots,w_5$ are contained in ${\rm Int}F$.
\end{enumerate}
Let $e=\Gamma_m\cap\partial F$.
Then $e$ is an internal edge of label $m$
connecting $w_1$ and $w_4$
(see Fig.~\ref{fig25}(a)).

{\bf Claim~3.}
$e\cap e_2'=\emptyset$
(see Fig.~\ref{fig25}(b)).

{\it Proof of Claim~$3$.}
Suppose 
$e\cap e_2'\not=\emptyset$.
By (1), the two white vertices $w_1$ and $w_4$ 
are in the outside of $E$. Thus 
the intersection $e\cap E$ consists of proper arcs of $E$.
Since $\partial E=e_2'\cup \widetilde e$ and
since $\widetilde e$ is an edge of label $m+1$,
the intersection $e\cap E$ consists of 
$(E,e_2')$-arcs of label $m$.
This contradicts Claim~2.
Thus $e\cap e_2'=\emptyset$.
Hence Claim~3 holds.
{\hfill {$\square$}\vspace{1.5em}}

By Claim~3,
the edge $e_2'$ is a proper arc of $F$.
Thus the edge $e_2'$ separates the region $F$
into the disk $E$ and a disk, say $E'$.
Moreover
\begin{enumerate}
\item[(4)] the boundary $\partial E'$ consists of
$e,e_2'$ and two internal edges of label $m+1$.
\end{enumerate}

Since the five white vertices $w_1,w_2,\cdots,w_5$
are contained in $\Gamma_{m+1}$
and since $\Gamma$ is of type $(m;2,3,2)$,
we have
\begin{enumerate}
\item[(5)] $w_6,w_7\in \Gamma_{m+2}\cap\Gamma_{m+3}$.
\end{enumerate}

{\bf Claim~4.}
We can assume $w(\Gamma\cap {\rm Int}E')=0$
by C-moves in a neighborhood of $E'$ 
keeping $\partial E'$ fixed.

{\it Proof of Claim~$4$.}
By (3),
none of $w_1,w_2,\cdots,w_5$
are contained in ${\rm Int}E'$.
Thus it is possible that
${\rm Int}E'$ contains $w_6$ or $w_7$.

Suppose that ${\rm Int}E'$ contains the white vertex $w_6$.
Since the edge $e_3'$ of label $m+2$ connects
$w_6$ and $w_3$,
the edge $e_3'$ intersects $\partial E'-w_3$.
Hence by (4),
the edge $e_3'$ intersects the edge $e$ of label $m$.
Let $x$ be the point in $e_3'$
with $e_3'[w_6,x]\cap e=x$.
Then by (5)
the vertex $w_6\in\Gamma_{m+2}\cap\Gamma_{m+3}$
 connects the point $x$ of the edge $e$ of label $m$
by the arc $e_3'[w_6,x]$ of label $m+2$.
Hence by Shifting Lemma (Lemma~\ref{Shift}),
we can shift the white vertex $w_6$ to the outside of $E'$
along the arc $e_3'[w_6,x]$ by
C-I-R2 moves, C-I-R3 moves and C-I-R4 moves
in a neighborhood of the arc $e_3'[w_6,x]$ 
keeping $\displaystyle{\bigcup_{i<0}\Gamma_{m+2+i}}$ fixed.
Thus we can shift the white vertex $w_6$ to 
the outside of $E'$ keeping $\partial E'$ fixed.

Similarly we can shift the white vertex $w_7$ 
to the outside of $E'$ keeping $\partial E'$ fixed.
Hence Claim~4 holds.
{\hfill {$\square$}\vspace{1.5em}}

By Claim~4,
we can apply New Disk Lemma 
(Lemma~\ref{NewDiskLemma}) for the disk $E'$
so that we have the following claim:

{\bf Claim~5.}
We can assume that
$\Gamma$ is $(E',e)$-arc free by C-moves 
in a neighborhood of the edge $e$ 
keeping $\partial E'$ fixed.

By the similar way of the proof of Claim~3,
we have the following claim:

{\bf Claim~6.}
$\Gamma_{m+2}\cap e=\emptyset$.

Now we are ready to prove Lemma~\ref{NoMinimalTypeH}.
Suppose that $\Gamma$ is minimal.
By Claim~6,
we can apply C-II moves along the edge $e$
for the black vertex in the terminal edge $e_1$
of label $m+1$ at $w_1$
such that 
the black vertex in $e_1$ is 
near the white vertex $w_4$
(see Fig.~\ref{fig25}(c)).
By a C-I-M2 move between $e_1$ and 
an internal edge of label $m+1$ 
connecting $w_4$ and $w_5$,
we obtain a new terminal edge of label $m+1$ at $w_4$
but not middle at $w_4$
(see Fig.~\ref{fig25}(d)).
This contradicts Assumption~\ref{AssumeTerminal}.
Therefore $\Gamma$ is not minimal.
We complete the proof of 
Lemma~\ref{NoMinimalTypeH}.
\end{Proof}



\section{Proof of Main Theorem}
\label{s:Main}

In this section,
we shall show the main theorem 
(Theorem~\ref{MainTheorem}).

\begin{lemma}
\LABEL{LemmaTypeHTypeHTypeA}
Let $\Gamma$ be a chart of type $(m;2,3,2)$.
If $\Gamma$ contains the pseudo chart as shown 
in Fig.~\ref{fig08}$($c$)$,
then $\Gamma$ is not minimal.
\end{lemma}

\begin{Proof}
Suppose that $\Gamma$ is minimal.
We use the notations as shown in 
Fig.~\ref{fig08}$($c$)$
where 
\begin{enumerate}
\item[(1)]
 $e_2,e_4$ are internal edges (possibly terminal edges)
 of label $m+2$ middle at $w_2,w_4$ and 
 oriented inward at $w_2,w_4$
 respectively.
\end{enumerate}
Since $w_1,w_3\in\Gamma_m\cap\Gamma_{m+1}$,
$w_2,w_4,w_5\in\Gamma_{m+1}$
and $\Gamma$ is of type $(m;2,3,2)$,
we have
\begin{enumerate}
\item[(2)]  $\Gamma_{m+1}\cap\Gamma_{m+2}=\{w_2,w_4,w_5\}$.
\end{enumerate}

By Corollary~\ref{TypeGTypeH} and 
Proposition~\ref{LemmaTypeGTypeB},
the graph $\Gamma_{m+2}$ contains the graph
as shown in Fig.~\ref{fig05}(b).
Hence by using Lemma~\ref{OriBWvertex},
the graph $\Gamma_{m+2}$ contains
one of the two graphs as shown in 
Fig.~\ref{fig26}(a),(b).
We also use the notations as shown in
Fig.~\ref{fig26}(a),(b)
where
$v_1,v_2,\cdots,v_5$ are the white vertices
in $\Gamma_{m+2}$, and
\begin{enumerate}
\item[(3)]
$\widetilde{e}_2,\widetilde{e}_4$ are internal edges 
 of label $m+2$ middle at $v_2,v_4$ respectively,
\item[(4)] the edge of label $m+2$ middle at $v_5$
is a terminal edge.
\end{enumerate}

{\bf Claim.}
$w_2\not= v_2$ and $w_2\not=v_4$.

{\it Proof of Claim.}
If $w_2=v_2$,
then $e_2$ is middle at $w_2(=v_2)$ by (1).
Thus by (3) we have $e_2=\widetilde{e}_2$.

Let $D$ be the 2-angled disk of $\Gamma_{m+1}$
with $\partial D\ni w_1,w_2$ and $D\supset e_2(=\widetilde{e}_2)$
(see Fig.~\ref{fig08}(c) 
and Fig.~\ref{fig26}(c)).
Since $w_1\in \Gamma_m\cap\partial D$,
the interior ${\rm Int}D$ contains $v_3,v_4$ and $v_5$
(see Fig.~\ref{fig26}(c)).
Thus $w(\Gamma\cap {\rm Int}D)\ge3$.
Since $\Gamma$ is of type $(m;2,3,2)$,
we have $w(\Gamma)=7$ and $w(\Gamma_{m+1})=5$.
Hence
$$7=w(\Gamma)\ge w(\Gamma_{m+1})+w(\Gamma\cap {\rm Int}D)\ge 5+3=8.$$
This is a contradiction.
Thus $w_2\not=v_2$.

Similarly we can show that
$w_2\not=v_4$.
Hence Claim holds. {\hfill {$\square$}\vspace{1.5em}}

Since $\Gamma_{m+2}$ contains the graph as shown
in Fig.~\ref{fig05}(b),
by Lemma~\ref{LemmaTypeGTypeHGammaM+2}
the chart $\Gamma$ contains one of 
the RO-families of the two pseudo charts 
as shown in Fig.~\ref{fig09}(c),(d).
Hence there are two cases:
(i) $v_1,v_3\in \Gamma_{m+3}$ or $v_3,v_5\in \Gamma_{m+3}$,
(ii) $v_1,v_4\in \Gamma_{m+3}$ or $v_2,v_5\in \Gamma_{m+3}$.

{\bf Case (i).}
Without loss of generality
we can assume $v_1,v_3\in \Gamma_{m+3}$.
Since $v_1,v_2,v_3,v_4,v_5\in\Gamma_{m+2}$ and
$\Gamma$ is of type $(m;2,3,2)$,
we have $\Gamma_{m+1}\cap\Gamma_{m+2}=\{v_2,v_4,v_5\}$. 
Thus by (2) we have
$\{v_2,v_4,v_5\}=\{w_2,w_4,w_5\}$.
Hence by Claim,
we have $w_2=v_5$
and 
\begin{enumerate}
\item[(5)] $w_4\in \{v_2,v_4\}$.
\end{enumerate}

Since the edge $e_2$ of label $m+2$ 
is middle at $w_2(=v_5)$ by (1),
the edge $e_2$
is the terminal edge at $v_5$
by (4).
Moreover by (1), 
the terminal edge $e_2$ is oriented inward at $w_2(=v_5)$.
Thus $\Gamma_{m+2}$ contains the graph as shown in
Fig.~\ref{fig26}(a).
Hence 
\begin{enumerate}
\item[(6)] the two edges $\widetilde{e}_2,\widetilde{e}_4$
 of label $m+2$
are oriented outward at $v_2,v_4$ respectively.
\end{enumerate}

On the other hand,
since the edge $e_4$ of label $m+2$
is middle at $w_4$ by (1),
either $e_4=\widetilde{e}_2$ or $e_4=\widetilde{e}_4$ 
by (3) and (5).
Since the edge $e_4$ of label $m+2$
is oriented inward at $w_4$ by (1),
 one of the two edges $\widetilde{e}_2,\widetilde{e}_4$ 
is oriented inward at $w_4$ $(v_2$ or $v_4)$.
This contradicts (6).
Hence Case (i) does not occur.

{\bf Case (ii).}
Without loss of generality
we can assume $v_1,v_4\in \Gamma_{m+3}$.
Since $v_1,v_2,v_3,v_4,v_5\in\Gamma_{m+2}$
and $\Gamma$ is of type $(m;2,3,2)$,
we have $\Gamma_{m+1}\cap\Gamma_{m+2}=\{v_2,v_3,v_5\}$.
Thus by (2)
we have 
\begin{enumerate}
\item[(7)] $\{v_2,v_3,v_5\}=\{w_2,w_4,w_5\}$.
\end{enumerate}

Let $e_5$ be the internal edge 
(possibly terminal edge)
of label $m+2$ middle at $w_5$,
and $\widetilde{e}_3,\widetilde{e}_5$ 
the terminal edges of 
label $m+2$ at $v_3,v_5$ respectively.
Since
$\widetilde{e}_3,\widetilde{e}_5$ 
are middle at $v_3,v_5$ respectively
by Assumption~\ref{AssumeTerminal} and 
since
$e_2,e_4,\widetilde{e}_2$ are middle at $w_2,w_4,v_2$ respectively by (1) and (3),
the condition (7) implies that 
\begin{enumerate}
\item[(8)] $\{\widetilde{e}_2,\widetilde{e}_3,\widetilde{e}_5\}=\{e_2,e_4,e_5\}$.
\end{enumerate}

Since there are two internal edges of label $m+1$
oriented from $w_5$ to $w_4$
(see Fig.~\ref{fig08}(c)),
the edge $e_5$ is oriented outward at $w_5$.
Moreover
since $e_2,e_4$ are oriented inward at $w_2,w_4$
respectively by (1),
the condition (8) implies that
the graph $\Gamma_{m+2}$
contains the graph as shown in 
Fig.~\ref{fig26}(b)
and $\{w_2,w_4\}=\{v_2,v_3\}$,
$w_5=v_5$.
Thus by Claim,
we have $w_2=v_3$ and $w_4=v_2$.
Hence the condition $w_4=v_2$
implies $e_4=\widetilde{e}_2$.
Thus
the edge $e_4(=\widetilde{e}_2)$ connects 
$w_4(=v_2)$ and $w_2(=v_3)$
(see Fig.~\ref{fig26}(b)).
However the edge $e_4$
must intersect the boundary 
of the 2-angled disk of $\Gamma_{m+1}$
(see Fig.~\ref{fig26}(d)).
This contradicts the definition of a crossing.
Thus Case (ii) does not occur.

Therefore $\Gamma$ is not minimal.
We complete the proof of 
Lemma~\ref{LemmaTypeHTypeHTypeA}.
\end{Proof}

\begin{figure}[htb]
\centerline{\includegraphics{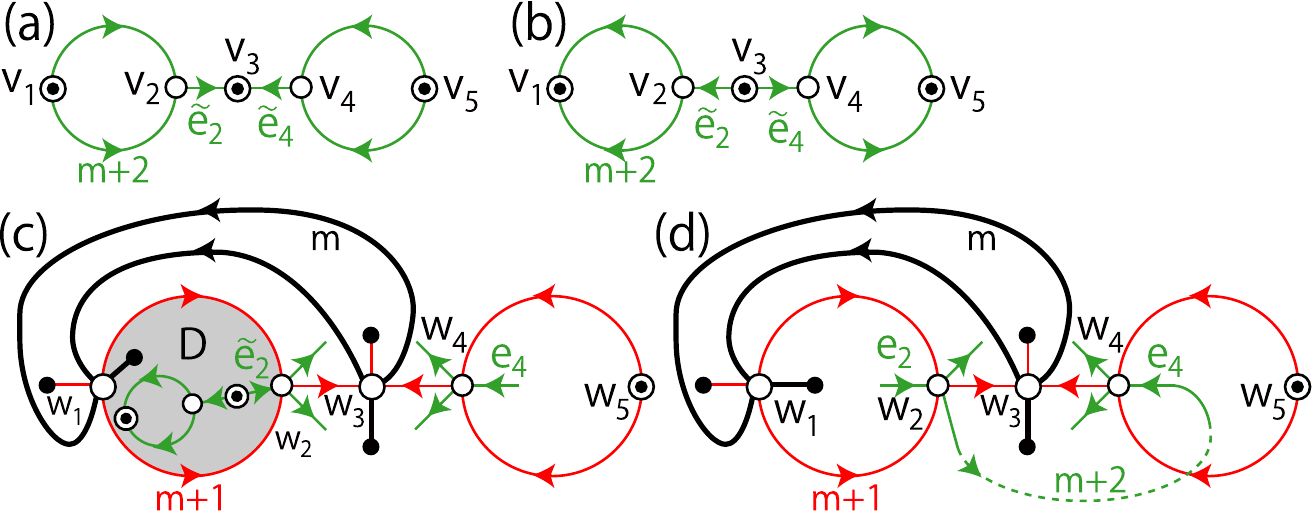}}
\caption{\LABEL{fig26}
(a), (b) The graphs as shown in Fig.~\ref{fig05}(b).
(c) The disk $D$ contains at least three white vertices 
in its interior.
(d) There is an internal edge of label $m+2$
connecting $w_2$ and $w_4$.
}
\end{figure}

{\it Proof of Theorem~\ref{MainTheorem}.}
Suppose that there exists a minimal chart
$\Gamma$ of type $(m;2,3,2)$.
By Corollary~\ref{TypeGTypeH} and
Proposition~\ref{LemmaTypeGTypeB}, 
each of $\Gamma_{m+1}$ and $\Gamma_{m+2}$
contains the graph as shown in 
Fig.~\ref{fig05}$($b$)$.
By Lemma~\ref{LemmaTypeH} and
Lemma~\ref{LemmaTypeHTypeHTypeA},
the chart $\Gamma$ contains one of 
the RO-family of the pseudo chart as shown
in Fig.~\ref{fig08}(d).
We use the notations as shown in 
Fig.~\ref{fig08}$($d$)$ where
\begin{enumerate}
\item[(1)]  
 $e_2$ is an internal edge (possibly terminal edge)
 of label $m+2$ middle at $w_2$ and 
 oriented inward at $w_2$.
\end{enumerate}
Since $w_1,w_4\in\Gamma_m\cap\Gamma_{m+1}$,
$w_2,w_3,w_5\in\Gamma_{m+1}$
and $\Gamma$ is of type $(m;2,3,2)$,
we have
\begin{enumerate}
\item[(2)]  $\Gamma_{m+1}\cap\Gamma_{m+2}=\{w_2,w_3,w_5\}$.
\end{enumerate}

Since $\Gamma_{m+2}$ contains the graph
as shown in Fig.~\ref{fig05}(b),
the graph $\Gamma_{m+2}$ contains 
one of the two graphs as shown in 
Fig.~\ref{fig26}(a),(b).
We also use the notations as shown in
Fig.~\ref{fig26}(a),(b),
where
$v_1,v_2,\cdots,v_5$ are the white vertices
in $\Gamma_{m+2}$,
\begin{enumerate}
\item[(3)]  $\widetilde{e}_2$ is middle at $v_2$.
\end{enumerate}

Moreover,
by Lemma~\ref{LemmaTypeGTypeHGammaM+2},
the chart $\Gamma$ contains one of the RO-families
of the two pseudo charts 
as shown in Fig.~\ref{fig09}(c),(d).
Hence there are two cases:
(i) $v_1,v_3\in \Gamma_{m+3}$ or $v_3,v_5\in \Gamma_{m+3}$,
(ii) $v_1,v_4\in \Gamma_{m+3}$ or $v_2,v_5\in \Gamma_{m+3}$.

{\bf Case (i).}
By changing the labels $\cdots,m,m+1,m+2,m+3,\cdots$ 
into $\cdots,m+3,m+2,m+1,m,\cdots$ respectively,
the chart $\Gamma$ contains the pseudo chart as shown
in Fig.~\ref{fig08}(c).
Thus by Lemma~\ref{LemmaTypeHTypeHTypeA},
we show that $\Gamma$ is not a minimal chart.
This is a contradiction.
Hence Case (i) does not occur.

{\bf Case (ii).}
Without loss of generality
we can assume $v_1,v_4\in \Gamma_{m+3}$.
Since $v_1,v_2,v_3,v_4,v_5\in\Gamma_{m+2}$ and
$\Gamma$ is of type $(m;2,3,2)$,
we have $\Gamma_{m+1}\cap\Gamma_{m+2}=\{v_2,v_3,v_5\}$.
Thus by (2),
we have 
\begin{enumerate}
\item[(4)] $\{v_2,v_3,v_5\}=\{w_2,w_3,w_5\}$.
\end{enumerate}

By the same way of the proof of Claim
in Lemma~\ref{LemmaTypeHTypeHTypeA},
we can show that 
\begin{enumerate}
\item[(5)] $w_2\not= v_2$ and $w_2\not=v_4$.
\end{enumerate}

Let $e_3,e_5$ be internal edges (possibly terminal edges)
of label $m+2$ middle at $w_3,w_5$ respectively.
Let $\widetilde{e}_3,\widetilde{e}_5$ be 
the terminal edges at $v_3,v_5$ respectively.
Since $\widetilde{e}_3,\widetilde{e}_5$ are middle 
at $v_3,v_5$ respectively
by Assumption~\ref{AssumeTerminal}
and since $e_2,\widetilde{e}_2$ are middle at $w_2,v_2$ 
respectively
by (1) and (3),
the condition (4) implies that
\begin{enumerate}
\item[(6)] $\{\widetilde{e}_2,\widetilde{e}_3,\widetilde{e}_5\}=\{e_2,e_3,e_5\}$.
\end{enumerate}

Since there are two internal edges of label $m+1$
oriented inward at $w_3$,
and since there are two internal edges of label $m+1$
oriented outward at $w_5$
(see Fig.~\ref{fig08}(d)),
\begin{enumerate}
\item[(7)] the edge $e_3$ is oriented inward at $w_3$
and the edge $e_5$ is oriented outward at $w_5$.
\end{enumerate}
Furthermore
since $e_2$ is oriented inward at $w_2$ by (1),
the conditions (6) and (7) imply 
that $\Gamma_{m+2}$ contains
the graph as shown in Fig.~\ref{fig26}(b).
Hence $\{w_2,w_3\}=\{v_2,v_3\}$,
$w_5=v_5$.
Thus
by (5),
we have $w_2=v_3$ and $w_3=v_2$.
Hence 
$e_2=\widetilde{e}_3$, $e_3=\widetilde{e}_2$ 
and $e_5=\widetilde{e}_5$.
Hence $e_3(=\widetilde{e}_2)\ni w_2(=v_3),w_3(=v_2)$, and
both of $e_2(=\widetilde{e}_3),e_5(=\widetilde{e}_5)$ 
are terminal edges.
Furthermore,
there exist two internal edges of label $m+2$
connecting $w_3(=v_2)$ and $v_1$,
and
there exist three internal edges of label $m+2$
at $v_4$
containing $w_2(=v_3),w_5(=v_5),w_5$ respectively.
Therefore
$\Gamma$ satisfies 
the condition of Lemma~\ref{NoMinimalTypeH}
(see Fig.~\ref{fig25}(a)).
By Lemma~\ref{NoMinimalTypeH},
the chart $\Gamma$ is not minimal.
This is a contradiction.
Hence Case (ii) does not occur.

Therefore both cases do not occur.
Hence there does not exist any minimal chart 
of type $(m;2,3,2)$.
We complete the proof of 
Theorem~\ref{MainTheorem}.
{\hfill {$\square$}\vspace{1.5em}}




\vspace{5mm}

\begin{minipage}{65mm}
{Teruo NAGASE
\\
{\small Tokai University \\
4-1-1 Kitakaname, Hiratuka \\
Kanagawa, 259-1292 Japan\\
\\
nagase@keyaki.cc.u-tokai.ac.jp
}}
\end{minipage}
\begin{minipage}{65mm}
{Akiko SHIMA 
\\
{\small Department of Mathematics, 
\\
Tokai University
\\
4-1-1 Kitakaname, Hiratuka \\
Kanagawa, 259-1292 Japan\\
shima@keyaki.cc.u-tokai.ac.jp
}}
\end{minipage}

\vspace{7mm}

{\bf List of terminologies}\vspace{2mm}\\
{\small $
\begin{array}{ll||}
\text{$k$-angled disk} & p8 \\
\text{BW-vertex} & p5 \\
\text{C-move equivalent} & p3 \\
\text{chart} & p2 \\
\text{complexity} & p3 \\
\text{$(D,\alpha)$-arc of label $k$} & p16 \\
\text{$(D,\alpha)$-arc free} & p16 \\
\text{feeler} & p8 \\
\text{free edge} & p3 \\
\text{hoop} & p4 \\
\text{internal edge} & p9 \\
\text{inward} & p3 \\
\text{inward arc} & p10 \\
\text{IO-Calculation} & p11 \\
\text{keeping $X$ fixed} & p12 \\
\text{local complexity} & p21 \\
\text{locally minimal} & p21 \\
\text{loop} & p23 \\
\end{array}
~~
\begin{array}{ll}
\text{middle arc} & p3 \\
\text{middle at $v$} & p3 \\
\text{minimal chart} & p3 \\
\text{outward} & p3 \\
\text{outward arc} & p10 \\
\text{oval} & p5 \\
\text{point at infinity $\infty$} & p4 \\
\text{proper arc} & p9 \\
\text{pseudo chart} & p6 \\
\text{ring} & p4 \\
\text{RO-family} & p6 \\
\text{simple hoop} & p4 \\
\text{skew $\theta$-curve} & p5 \\
\text{special $k$-angled disk} & p21 \\
\text{special oval} & p12 \\
\text{terminal edge} & p4 \\
\text{type $(m;n_1,n_2,\cdots,n_k)$} & p1 \\
\text{$w$ connects with $p$ by an arc $\beta$} & p28 \\
\end{array}
$}

\newpage

\vspace{0.5cm}

{\bf List of notations}\vspace{2mm}\\
{\small $
\begin{array}{ll}
\text{$\Gamma_m$} & p1 \\
\text{${\rm Int}X$} & p5 \\
\text{$\partial X$} & p5 \\
\text{$Cl(X)$} & p5 \\
\text{$\partial \alpha$} & p5 \\
\text{${\rm Int}\alpha$} & p5 \\
\text{$w(X)$} & p5 \\
\text{$a_{ij}$,$b_{ij}$} & p6 \\
\text{$c(X)$} & p20 \\
\text{$\ell c(D;\Gamma)$} & p21\\
\text{$\alpha[p,q]$} & p28 \\
\end{array}
$
}


\begin{thebibliography}{99}


\bibitem{KnottedSurfaces}J.\ S.\ Carter and M.\ Saito, {"Knotted surfaces and their diagrams"}, Mathematical Surveys and Monographs, 55, American Mathematical Society, Providence, RI, (1998). MR1487374 (98m:57027)


\bibitem{H1} \textsc{I.\ Hasegawa}, 
{\it The lower bound of the w-indices of non-ribbon surface-links}, 
Osaka J. Math. {\bf 41} (2004), 891--909.
MR2116344 (2005k:57045)

\bibitem{INS} S.\ Ishida, T.\ Nagase and A.\ Shima, 
{\it Minimal $n$-charts with four white vertices}, J. Knot Theory Ramifications {\bf 20}, 689--711 (2011).
MR2806339 (2012e:57044)



\bibitem{BraidThree} S.\ Kamada, {\it Surfaces in $R\sp 4$ of braid index three are ribbon}, J. Knot Theory Ramifications {\bf 1} No. 2 (1992), 137--160. 
MR1164113 (93h:57039)

\bibitem{BraidBook} S.\ Kamada, {"Braid and Knot Theory in Dimension Four"}, 
Mathematical Surveys and Monographs,
Vol. 95, American Mathematical Society, (2002). MR1900979 (2003d:57050)


\bibitem{ChartApp1} T.\ Nagase and A.\ Shima, 
{\it Properties of minimal charts and their applications I}, J. Math. Sci. Univ. Tokyo {\bf 14} (2007), 69--97.
MR2320385 (2008c:57040)

\bibitem{ChartAppII} T.\ Nagase and A.\ Shima, 
{\it Properties of minimal charts and their applications II}, Hiroshima Math. J. {\bf 39} (2009), 1--35.
MR2499196 (2009k:57040)

\bibitem{ChartAppIII} T.\ Nagase and A.\ Shima, 
{\it Properties of minimal charts and their applications III}, Tokyo J. Math. {\bf 33} (2010), 373--392.
MR2779264 (2012a:57033)


\bibitem{ChartAppIV} T.\ Nagase and A.\ Shima, 
{\it Properties of minimal charts and their applications IV: Loops}, 
 J. Math. Sci. Univ. Tokyo {\bf 24} (2017), 
 195--237,
arXiv:1603.04639. MR3674447


\bibitem{ChartAppV} T.\ Nagase and A.\ Shima, 
{\it Properties of minimal charts and their applications V: charts of type $(3,2,2)$}, J. Knot Theory Ramifications
{\bf 28} No. 14 (2019) 1950084 (27 pages),
arXiv:1902.00007.

\bibitem{ChartAppVI} T.\ Nagase and A.\ Shima, 
{\it Properties of minimal charts and their applications VI: the graph $\Gamma_{m+1}$ in a chart $\Gamma$ of type $(m;2,3,2)$}, preprint, arXiv:2003.11909.

\bibitem{ChartAppVIII} T.\ Nagase and A.\ Shima, 
{\it Properties of minimal charts and their applications VIII-}, in preparation.

\bibitem{Chart33} T.\ Nagase and A.\ Shima, 
{\it Minimal charts of type $(3,3)$}, 
Proc. Sch. Sci. TOKAI UNIV. {\bf 52} (2017), 1--25,
 arXiv:1609.08257v2.

\bibitem{StI} T.\ Nagase and A.\ Shima, 
{\it The structure of a minimal $n$-chart with two crossings I: Complementary domains of $\Gamma_1\cup \Gamma_{n-1}$}, 
J. Knot Theory Ramifications {\bf 27} No. 14 (2018) 1850078 (37 pages), arXiv:1704.01232v3.

\bibitem{NST} T.\ Nagase, A.\ Shima and H.\ Tsuji, 
{\it The closures of surface braids obtained from minimal $n$-charts with four white vertices}, 
J. Knot Theory Ramifications {\bf 22} No. 2 (2013) 1350007 (27 pages). MR3037298 

\bibitem{ONS}
{M.\ Ochiai, T.\ Nagase and A.\ Shima}, 
{\it There exists no minimal $n$-chart with five white vertices}, 
Proc. Sch. Sci. TOKAI UNIV. 40 (2005), 1--18. 
MR2138333 (2006b:57035)

\bibitem{Tanaka} K.\ Tanaka, 
{\it A Note on CI-moves}, 
Intelligence of low dimensional topology 2006, 307--314, 
Ser. Knots Everything, 40, World Sci. Publ., Hackensack, NJ, 2007. MR2371740 (2009a:57017)
\end{thebibliography}
\end{document}